\documentclass[11pt,a4paper,twoside]{article}

\usepackage{a4}

\usepackage{amsmath}
\usepackage{amssymb}   
\usepackage{amsthm}

\usepackage[all]{xy}
\CompileMatrices
\UseComputerModernTips

\xyoption{ps}
\xyoption{dvips}

\newcommand{\Res}{\text{Res}}
\newcommand{\Gal}{\text{Gal}}
\newcommand{\GL}{\text{GL}}
\newcommand{\SL}{\text{SL}}
\newcommand{\PSL}{\text{PSL}}

\newcommand{\Eis}{\text{Eis}}
\newcommand{\Hom}{\text{Hom}}
\newcommand{\Ext}{\text{Ext}}
\newcommand{\CHE}{\text{CHE}}
\newcommand{\cusp}{\text{cusp}}
\newcommand{\Spec}{\text{Spec}}
\newcommand{\Pic}{\text{Pic}}
\newcommand{\CH}{\text{CH}}

\theoremstyle{plain}
\newtheorem{Thm}{Theorem}[section]

\newtheorem{Lem}[Thm]{Lemma}

\newtheorem{Cor}[Thm]{Corollary}

\newtheorem{Prop}[Thm]{Proposition}

\theoremstyle{definition}
\newtheorem{Rmk}[Thm]{Remark}

\newtheorem{Def}{Definition}

\theoremstyle{plain}
\newtheorem{T}{Theorem}

\theoremstyle{plain}

\begin{document}


\title{Realisations of Kummer-Chern-Eisenstein-Motives}
{\scriptsize}\author{Alexander Caspar}

\maketitle



\begin{abstract}
Inspired by the work of G. Harder (\cite{HaICM}, \cite{HaLNM}, \cite{HaMM}) we
construct via the motive of a Hilbert modular surface an extension
of a Tate motive by a Dirichlet motive. We compute the realisation
classes and indicate how this is linked to the Hodge-$1$-motive of the Hilbert
modular surface. 

\emph{Classification:} 11F41, 11F80, 14G35, 11G18
\end{abstract}

\section*{Introduction}

Let $D$ be a prime with $D \equiv 1
 (\text{mod} \: 4)$
 and let $F = \mathbb{Q}(\sqrt{D})$ be the real quadratic number field of
 discriminant $D$. 
 We choose $\sqrt{D} > 0$ and consider $F$ as a subfield of $\mathbb{R}$. We
 assume that the class number in the narrow sense $h^+$ is $1$. Denote by $\mathcal{O}_F \subset F$ the ring of integers. Its group of units $\mathcal{O}_F^*$ is isomorphic to $\mathbb{Z}
 \times (\mathbb{Z}/2\mathbb{Z})$ and there is a fundamental unit $\varepsilon_0\in\mathcal{O}_F^*$
 with norm 
$-1$, where we choose $\varepsilon_0 >1$. Moreover, we define
$\varepsilon:=\varepsilon_0^2$, which is a generator of the totally positive units. 
By $\chi_D(\cdot) := \left(\frac{D}{\cdot}\right)$ we
 denote the primitive Dirichlet
 character mod $D$ and $\zeta_F(-1)$ is the Dedekind $\zeta$-function of $F$
 at $-1$.

Consider the  Hilbert modular surface $S/\mathbb{Q}$ with full level.
Let 
$j:S\hookrightarrow\widetilde{S}$ be the open embedding of $S/\mathbb{Q}$ into the smooth toroidal compactification $\widetilde{S}/\mathbb{Q}$.
Then this gives rise to an exact sequence of motives
$ 0 \rightarrow \mathbb{Q}(0) \rightarrow H^2_c(S,\mathbb{Q}) \rightarrow H^2(S,\mathbb{Q}) \rightarrow \mathbb{Q}(-2) \rightarrow 0,$
where $\mathbb{Q}(n)$ is the Tate motive of weight $-2n$. 
With the interior cohomology
$H^2_!(S,\mathbb{Q}):=\text{Im}(H^2_c(S,\mathbb{Q})\rightarrow H^2(S,\mathbb{Q}))$ we get an extension 
$0\rightarrow H^2_!(S,\mathbb{Q})\rightarrow H^2(S,\mathbb{Q})\rightarrow \mathbb{Q}(-2)\rightarrow 0.$
Now the bottom $H^2_!(S,\mathbb{Q})$ is a semi-simple module for the Hecke
algebra and we can find a direct summand $H^2_{\CH}(S,\mathbb{Q}(1)) =
\mathbb{Q}(0) \oplus \mathbb{Q}(0)\chi_D,$ where $\mathbb{Q}(0)\chi_D$ is the Dirichlet motive  for the quadratic character
$\chi_D$ in the sense of Deligne. 
It is generated by
the first Chern class $c_1(L_1
\otimes L_2^{-1})$, where $L_1
\otimes L_2$ is the line bundle of Hilbert modular forms of weight
$(2,2)$. Hence one is faced with 
$0 \rightarrow \mathbb{Q}(0)\chi_D \rightarrow H^2_{\CHE}(S,\mathbb{Q}(1)) \rightarrow \mathbb{Q}(-1) \rightarrow 0,$
i.e. an element 
$[H^2_{\CHE}(S,\mathbb{Q}(1))(-)] \in \Ext^1_{\mathcal{MM}_{\mathbb{Q}}}(\mathbb{Q}(-1),\mathbb{Q}(0)\chi_D),$
which we call a \emph{Kummer-Chern-Eisenstein motive}.

Now we have the various realisations of our motive. There are the $l$-adic realisations
$[H^2_{\CHE,l}(S,\mathbb{Q}(1))(-)] \in
\Ext^1_{\mathcal{MGAL}}(\mathbb{Q}_l(-1),\mathbb{Q}_l(0)\chi_D)$, which are
mixed $\Gal(\overline{\mathbb{Q}}/\mathbb{Q})$-modules and we know that the latter group is the subgroup of norm one
elements in  $l$-adic completion of $F^*$ tensorised with $\mathbb{Q}_l$. Let $\widetilde{\varepsilon}:= \varepsilon^{-\frac{1}{2}
  \zeta_F(-1)^{-1}}$. Then we have 
\begin{T}
Let 
$
[H^2_{\CHE,l}(S,\mathbb{Q}(1))(-)]$
be the $l$-adic realisations of our Kummer-Chern-Eisenstein motive. Then $[H^2_{\CHE,l}(S,\mathbb{Q}(1))(-)]$ is $\widetilde{\varepsilon}.$ The corresponding $l$-adic Galois representation comes from the Kummer field extension $F\left(\sqrt[l^{\infty}]{\widetilde{\varepsilon}}, \zeta_{l^{\infty}}\right)$ of $F$ attached to $\widetilde{\varepsilon}$, i.e.
$\Gal(\overline{\mathbb{Q}}/\mathbb{Q})\ni\sigma \mapsto
\left(\begin{smallmatrix} \chi_D(\sigma) & \tau_{\widetilde{\varepsilon}}(\sigma) \alpha^{-1}(\sigma) \\
0 & \alpha^{-1}(\sigma) \end{smallmatrix}\right),$ 
where $\alpha$ is the cyclotomic character and $\tau_{\widetilde{\varepsilon}}(\sigma)$ is defined by 
$
\frac{\sigma\left(\sqrt[l^{\infty}]{\widetilde{\varepsilon}}\right)}{\sqrt[l^{\infty}]{\widetilde{\varepsilon}}}
= \zeta_{l^{\infty}}^{\tau_{\widetilde{\varepsilon}}(\sigma)}.$ 
\end{T}
There is also the Hodge-de Rham extension, which is a mixed Hodge-de Rham structure
$[H^2_{\CHE,\infty}(S,\mathbb{Q}(1))(-)] \in
\Ext^1_{\mathcal{MH}d\mathcal{R}_{\mathbb{Q}}}(\mathbb{Q}(-1),\mathbb{Q}(0)\chi_D)$ 
and we can identify the latter group with $\mathbb{R}$. 
\begin{T}
Let $[H^2_{\CHE,\infty}(S,\mathbb{Q}(1))(-)]$
be the Hodge-de Rham realisation of our motive. Then $[H^2_{\CHE,\infty}(S,\mathbb{Q}(1))(-)]$ is 
$\log \widetilde{\varepsilon}.$ 
\end{T} 
These theorems tell us that the realisations of our motive are exactly those
of the Kummer motive $K \langle \widetilde{\varepsilon} \rangle$ attached to $\widetilde{\varepsilon}$,
$$\left([H^2_{\CHE,\infty}(S,\mathbb{Q}(1))(-)],[H^2_{\CHE,l}(S,\mathbb{Q}(1))(-)]\right) = (\log \widetilde{\varepsilon}, \widetilde{\varepsilon}) = (K \langle \widetilde{\varepsilon} \rangle_{\infty},K \langle \widetilde{\varepsilon} \rangle_l),$$
but we  do not know, whether actually $[H^2_{\CHE}(S,\mathbb{Q}(1))(-)] =
K \langle \widetilde{\varepsilon} \rangle$. Such a Kummer motive $K \langle \widetilde{\varepsilon} \rangle$ is isomorphic
to a one-motive $M_{\widetilde{\varepsilon}}$ attached to
$\widetilde{\varepsilon}$. 
We can appoint a good candidate for $M_{\widetilde{\varepsilon}}$. Consider
$\widetilde{L}=\widetilde{L_1}^{-1} \otimes \widetilde{L_2}$, where
$\widetilde{L}$ is the unique prolongation of  $L_i$ to the compact surface
$\widetilde{S}$, which has trivial Chern class on the boundary, that is, if we
denote by $u: \text{Pic}(\widetilde{S}) \rightarrow
\text{Pic}(\widetilde{S}_{\infty})$ the restriction map, we have $u(\widetilde{L_1}^{-1} \otimes \widetilde{L_2}) \in \text{Pic}^0(\widetilde{S}_{\infty})$. We know that
$\text{Pic}^0(\widetilde{S}_{\infty}) \simeq \mathbb{G}_m$ and we get that 
$[\mathbb{Z}(\chi_D) \cdot \widetilde{L} \stackrel{u}{\rightarrow} \text{Pic}^0(\widetilde{S}_{\infty})]$
is a $1$-motive, where $\mathbb{Z}(\chi_D)
\cdot \widetilde{L}$ is the submodule generated by  $\widetilde{L}$. 
\begin{T}
Consider the motive $T([\mathbb{Z}(\chi_D) \cdot \widetilde{L}
\stackrel{u}{\rightarrow} \text{Pic}^0(\widetilde{S}_{\infty})]^{\vee})
\otimes \mathbb{Q}$ of the dual Kummer-$1$-motive $[\mathbb{Z}(\chi_D)
\cdot \widetilde{L} \stackrel{u}{\rightarrow}
\text{Pic}^0(\widetilde{S}_{\infty})]^{\vee}.$ 
Then $[H^2_{\CHE,l}(S,\mathbb{Q}(1))(-)] \simeq T_l([\mathbb{Z}(\chi_D) \cdot
\widetilde{L} \stackrel{u}{\rightarrow}
\text{Pic}^0(\widetilde{S}_{\infty})]^{\vee}) \otimes \mathbb{Q}_l.$ 
In particular, we have
$[H^2_{\CHE,l}(S,\mathbb{Q}(1))(-)] \simeq T_l(M_{\widetilde{\varepsilon}}) \otimes \mathbb{Q}_l.$
\end{T}
%
There is a third one-motive, the Hodge-one-motive $\eta_S$ attached
to the surface $S$. It comes from the Hodge structure given by the second cohomology of $S$. Using \cite{Ca} and the previous theorem we obtain.
\begin{T}
The Kummer-$1$-motive attached to $\varepsilon^{-2}$ is isomorphic to a submotive of the Hodge-one-motive $\eta_{S}$. In particular, 
the dual of the Kummer-Chern-Eisenstein motive $[H^2_{\CHE}(S,\mathbb{Q}(1))(-)]^{\vee}$ is isomorphic to a submotive of the realisation $T(\eta_S) \otimes \mathbb{Q}$ of the Hodge-one-motive.
\end{T} 
\emph{Acknowledgement.} These are improved parts of the author's thesis, which was
written at the MPI f\"ur Mathematik in Bonn under the supervision of
Prof. G. Harder. I am indebted to Prof. Harder for the guidance and advise in mathematics
from the very beginning. Moreover, I would like to thank the MPIM in Bonn for the extraordinary
support and working conditions. 

\section{The Kummer-Chern-Eisenstein Motive}
Let $F$ be the real quadratic number field with the conventions and notations
that we fixed in the beginning of the introduction. Denote by $\Theta \in \Gal(F/\mathbb{Q})$ the nontrivial
element in the Galois group of $F/\mathbb{Q}$, and write $x^{\Theta}$ for the
Galois conjugate of $x \in F$. 
\subsection{Hilbert Modular Surfaces}\label{HMS}
We start with the algebraic group 
$G/\mathbb{Q}:=\Res_{F/\mathbb{Q}}(\GL_2/F),$ 
i.e. we have $G_{\infty}:=G(\mathbb{R}) \simeq \GL_2(\mathbb{R})\times \GL_2(\mathbb{R}),$ 
and $G(\mathbb{Q})=\GL_2(F)$. Note that  $\Theta$ interchanges the two copies of $\GL_2(\mathbb{R})$.
Define 
$K_{\infty}:=\left\{\left(\begin{smallmatrix} a_1 & -b_1 \\ b_1 & a_1 \end{smallmatrix}\right) \times \left(\begin{smallmatrix} a_2 & -b_2 \\ b_2 & a_2 \end{smallmatrix}\right)\right\} \subset G^0_{\infty},$
where $G^0_{\infty} \subset  G_{\infty}$ is the connected component of the identity.
The quotient is 
$X:=G_{\infty}/K_{\infty} \simeq (\mathfrak{H}^+ \cup \mathfrak{H}^-) \times (\mathfrak{H}^+ \cup \mathfrak{H}^-),$ 
where $\mathfrak{H}^+ \subset \mathbb{C}$, resp. $\mathfrak{H}^- \subset
\mathbb{C}$ is the upper half plane, resp. the lower half plane. We write in
the following  $\mathfrak{H} = \mathfrak{H}^+$. Let $\mathbb{A}_f=\widehat{\mathbb{Z}}\otimes\mathbb{Q}$ be the ring of finite adeles of $\mathbb{Q}$.
For an open compact subgroup $K_f \subset G(\mathbb{A}_f)$ of the finite adelic points of $G$ we define 
$S_{K_f}(\mathbb{C}):=G(\mathbb{Q})\backslash \left(X \times G(\mathbb{A}_f)/K_f \right),$
and we get a quasi-projective complex algebraic surface, which one calls a
\emph{Hilbert modular surface}.

We choose the standard maximal compact subgroup 
$K_f=K_0=\prod_{\mathfrak{p}}\GL_2(\mathcal{O}_{\mathfrak{p}}).$ This has the
effect that $S(\mathbb{C}):=S_{K_0}(\mathbb{C})$ is connected and moreover,
$S_{K_0}(\mathbb{C})=  {\Gamma} \backslash (\mathfrak{H} \times
\mathfrak{H}),$ with $\Gamma := \PSL_2(\mathcal{O}_F).$ 
Now $S_{K_0}(\mathbb{C})$ has cyclic (quotient-) singularities, caused by the
fixed points of the $K_0$-action. These are mild singularities, that can be
resolved by finite chains of rational curves, cf. \cite{vdG}, II.6. Let us denote this resolution again by $S_{K_0}(\mathbb{C})$.
To see that $S_{K_f}(\mathbb{C})$ are  the complex points of a scheme $S_{K_f}$ defined over $\mathbb{Q}$, one has to invoke the theory of canonical models by Deligne and Shimura (see e.g. \cite{De-Shimura}). 
This gives that $S_{K_f}$ is a quasi-projective scheme over $\mathbb{Q}$. 
Recall that this scheme $S/\mathbb{Q}$ represents the coarse moduli space of
polarised abelian surfaces with real multiplication by $\mathcal{O}_F$. To get
a fine moduli space (i.e. a universal family of abelian surfaces), one has to
introduce a level $N$-structure. Otherwise one is just left with a moduli
stack $\mathcal{S}$. For this we consider the congruence subgroup $K_N := \{ g \in K_0 | g \equiv \text{id} \: \text{mod} \: N \}.$ 
Now according to \cite{Ra}, the scheme $S_{K_N}/\mathbb{Q}$  represents for $N
\geq 3$ the fine moduli space of polarised abelian surfaces with real
multiplication by $\mathcal{O}_F$ and level $N$-structure.

There are different ways to compactify our surface $S$. First, there is the (singular) Baily-Borel compactification $S \hookrightarrow \overline{S}$. Basically this is
$\overline{S} = S\cup (\PSL_2(\mathcal{O}_F)\backslash \mathbb{P}_F^1),$ 
and here in our case (with class number $h = 1$ and full level $K_0$) it is
just given by adding a singular point $\infty$ at infinity (a cusp). 
This cusp singularity can be resolved in
a canonical way $\widetilde{S}\rightarrow \overline{S}$. For the complex
surface $S(\mathbb{C})$ this is due to Hirzebruch (as in
\cite{HiHMS}), but it can also be done over the rational numbers $\mathbb{Q}$
according to Rapoport (\cite{Ra} - see also \cite{HLR}). This resolution
$\widetilde{S}$ is the smooth toroidal compactification (\cite{AMRT}). The
boundary $\widetilde{S}_{\infty} :=\widetilde{S} - S$ is a polygon with
rational components $\widetilde{S}_{\infty,i} \simeq \mathbb{P}^1$, $i = 0,
..., n-1$. 
\subsection{The Line Bundles of Modular Forms} \label{clb}%
%
Let us start with the complex situation (i.e. with the Betti realisation). By
\cite{vdG}, II.7, we know that the line bundle of modular forms is a product
$L_1 \otimes L_2$, where each factor corresponds to the factor of automorphy
$(c z_1 + d)^2$, resp. $(c^{\Theta} z_2 + d^{\Theta})^2$. The sections of $L_1 \otimes L_2$ on $S(\mathbb{C})$ are Hilbert modular forms
$f(z_1,z_2)$ of weight $(2,2)$, and we know that $L_1 \otimes L_2 =
\Omega^2_{S(\mathbb{C})}$. We have the Chern class map
$c_1 : \Pic(S(\mathbb{C})) \rightarrow H^2(S(\mathbb{C}),\mathbb{Q}(1)),$
where $\mathbb{Q}(1) = 2\pi i \:\mathbb{Q}$. The map $c_1$ is induced by the exponential sequence
 $0 \rightarrow \mathbb{Z}(1) \rightarrow \mathcal{O}_{S(\mathbb{C})} \stackrel{\text{exp}}{\rightarrow}  \mathcal{O}^*_{S(\mathbb{C})} \rightarrow 0$
of sheaves on $S(\mathbb{C}),$ where $\Pic(S(\mathbb{C})) =
H^1(S(\mathbb{C}),\mathcal{O}_{S(\mathbb{C})}^*)$. Now we want to extend the
line bundles $L_1$, $L_2$ on $S(\mathbb{C})$ to the compact surface
$\widetilde{S}(\mathbb{C})$. This is given by
\begin{Lem}\label{LEMLIFT}

Let $L_i \in \Pic(S(\mathbb{C}))$,  $i = 1,2$ be the line bundles on
$S(\mathbb{C})$ as above. Then there is a unique line bundle $\widetilde{L_i}$ in
$\Pic(\widetilde{S}(\mathbb{C}))$ with trivial Chern class on the boundary, $\widetilde{L_i}|_{\widetilde{S}_{\infty}(\mathbb{C})} \in \text{Pic}^0(\widetilde{S}_{\infty}(\mathbb{C})),$ 
such that its restriction to the open part $S(\mathbb{C})$ is $L_i$. In
particular, with the interior cohomology $H^2_!(S(\mathbb{C}),\mathbb{Q}(1)) := \text{Im} \left(
  H^2_c(S({\mathbb{C}}),\mathbb{Q}(1)) \rightarrow
  H^2(S(\mathbb{C}),\mathbb{Q}(1)) \right)$ we have $c_1(L_i) \in H^2_!(S(\mathbb{C}),\mathbb{Q}(1)).$
\end{Lem} 
\begin{proof}
For the construction of the extensions $\widetilde{L_1}$, $\widetilde{L_2}$
compare e.g. \cite{vdG}, IV.2.
The tensor product is 
$\widetilde{L_1} \otimes \widetilde{L_2} =
\Omega^2(\log\widetilde{S}_{\infty}(\mathbb{C})),$ 
where $\Omega^2(\log\widetilde{S}_{\infty}(\mathbb{C}))$ denotes the sheaf of
differentials on $\widetilde{S}(\mathbb{C})$, which may have poles of at most
simple order along the boundary $\widetilde{S}_{\infty}(\mathbb{C})$. To proof
the uniqueness, we observe that there is, according to \cite{HLR}, Lemma 2.2, the following commutative triangle 
$$\xymatrix{
H^2_c(S(\mathbb{C}),\mathbb{Q}(1)) \ar[d] \ar[dr] &  \\
H^2(\widetilde{S}(\mathbb{C}),\mathbb{Q}(1)) \ar[r] & H^2(S(\mathbb{C}),\mathbb{Q}(1)),}$$
i.e. 
$H^2_!(S(\mathbb{C}),\mathbb{Q}(1)) =
\text{Im}\left(H^2(\widetilde{S}(\mathbb{C}),\mathbb{Q}(1)) \rightarrow
  H^2(S(\mathbb{C}),\mathbb{Q}(1)) \right).$ 
Now the Chern class $c_1 (\widetilde{L_i}) \in H^2(\widetilde{S}(\mathbb{C}),\mathbb{Q}(1))$ of the above line bundle is a preimage of $c_1(L_i)$, that has trivial Chern class on the boundary, i.e.
 $\widetilde{L_i}|_{\widetilde{S}_{\infty}(\mathbb{C})} \in
 \text{Pic}^0(\widetilde{S}_{\infty}(\mathbb{C})).$ 
We have to see that this extension is indeed unique. We prove that the extension cannot be modified by a divisor, that is supported on the boundary $\widetilde{S}_{\infty}(\mathbb{C})$.
For this we use the fact that the intersection matrix $(S_i \cdot S_j)_{ij}$
 of the boundary divisor $\widetilde{S}_{\infty}(\mathbb{C})$ is negative
 definite. In particular, the self-intersection number $S_i^2$ is at least
 $-2$, see for example \cite{vdG}, II.3. 
This means that a divisor, which is only supported on
$\widetilde{S}_{\infty}(\mathbb{C})$ has non zero degree. And hence we cannot
modify a line bundle with trivial Chern class on
$\widetilde{S}_{\infty}(\mathbb{C})$, for example our above $\widetilde{L_i}$,
 by a boundary divisor. Compare also \cite{HLR}, Hilfssatz 2.3.
\end{proof}
The sections of $\widetilde{L_1} \otimes \widetilde{L_2}$ over
$\widetilde{S}(\mathbb{C})$ are now \emph{meromorphic} Hilbert modular forms
of weight $(2,2)$. The product $\widetilde{L_1} \otimes \widetilde{L_2}$ is
trivial on the polygon at infinity, since there is, up to a constant, a
non-vanishing section, see \cite{vdG}, III. Lemma 3.2. But the restriction of each factor $\widetilde{L_i}$ to $\widetilde{S}_{\infty}(\mathbb{C})$ is not trivial. We see this in
\begin{Lem}\label{LEMCLB}
 
Let $ \widetilde{L_1} \otimes  \widetilde{L_2}$ be the line bundle of
meromorphic Hilbert modular forms of weight $(2,2)$ on $\widetilde{S}
(\mathbb{C})$. Let  $\varepsilon = \varepsilon_0^2 \in \mathcal{O}_F^*$ be as
above. Then the restriction
$\widetilde{L_i}|_{\widetilde{S}_{\infty}(\mathbb{C})} \in
\Pic^0(\widetilde{S}_{\infty}(\mathbb{C})) \simeq \mathbb{C}^*$ of each
factor $\widetilde{L_i}$ to the boundary $\widetilde{S}_{\infty}(\mathbb{C})$
is $
\widetilde{L_1}|_{\widetilde{S}_{\infty}(\mathbb{C})}=\varepsilon$ and $\widetilde{L_2}|_{\widetilde{S}_{\infty}(\mathbb{C})}=\varepsilon^{-1}.$
\end{Lem}
\begin{proof}
Restrict $\widetilde{L_i}$ to the boundary $\widetilde{S}_{\infty}(\mathbb{C})$
and use the explicit glueing of the components  
$\widetilde{S}_{\infty,i}(\mathbb{C})$ (cf. \cite{HiHMS}, 2.3) to see what
happenend to a section of
$\widetilde{L_i}|_{\widetilde{S}_{\infty}(\mathbb{C})}$ by going around the
polygon $\widetilde{S}_{\infty}(\mathbb{C})$. Moreover, we use that fact that
the units identify $\widetilde{S}_{\infty,i}(\mathbb{C})$ and
$\widetilde{S}_{\infty,i+n}(\mathbb{C})$, see loc. cit. or \cite{AMRT}, I.5.
\end{proof}
Note that the exponent $\pm1$ of $\varepsilon$ has been fixed by the orientation of
the boundary. 

To discuss the other realisations, we must
look at the moduli interpretation. For the fine moduli space there is the
universal family  $\mathcal{A}/ S_{K_N}$ of abelian surfaces, as above, with the zero section $s : S_{K_N} \rightarrow \mathcal{A}$. 
Let $\Omega^1_{\mathcal{A}/S_{K_N}}$ be the sheaf of relative differentials. Then we have the Lie algebra $\text{Lie}(\mathcal{A}) := s^*\Omega^1_{\mathcal{A}/S_{K_N}}$, which is a locally free $\mathcal{O}_F \otimes \mathcal{O}_{S_{K_N}}$-module of rank one (here we denote by $\mathcal{O}_{S_{K_N}}$ the structure sheaf of $S_{K_N}$).
And we define, as in \cite{Ra}, 6, a line bundle $\omega_{S_{K_N}}$ on $S_{K_N}$ by
$\omega_{S_{K_N}} := \text{Nm}_{\mathcal{O}_F \otimes
  \mathcal{O}_{S_{K_N}}/\mathcal{O}_{S_{K_N}}}\left(
  \text{Lie}(\mathcal{A})^{\vee} \right),$ where $\text{Nm}_{\mathcal{O}_F
  \otimes \mathcal{O}_{S_{K_N}}/\mathcal{O}_{S_{K_N}}}$ is the norm map from
$\mathcal{O}_F \otimes \mathcal{O}_{S_{K_N}}$ to $ \mathcal{O}_{S_{K_N}}$ and
$\text{Lie}(\mathcal{A})^{\vee}$ the dual Lie algebra. For its global sections we have 
\begin{Lem}

Let $L_1 \otimes L_2$ be the line bundle of Hilbert modular forms of weight $(2,2)$ on $S_{K_N}(\mathbb{C})$. Then its global sections are global sections of $\omega_{S_{K_N}(\mathbb{C})}^{\otimes 2}$ on $S_{K_N}(\mathbb{C})$, i.e. the line bundle $L_1 \otimes L_2$ comes from a line bundle on $S_{K_N}$.
\end{Lem}
\begin{proof}
This is \cite{Ra}, Lemme 6.12. Note that this is true for any level $N$, so even for the stack.
\end{proof} 
Still we have to show that each factor $L_i$ is defined over $F$, i.e. a line
bundle on $S_{K_N} \times F$. For this we base change our scheme $S_{K_N}$ over $\mathbb{Q}$ to $F$ and observe that $\mathcal{O}_F \otimes F$ decomposes into $F \oplus F$ along the action of $\Gal(F/\mathbb{Q}) = \{\text{id}, \Theta\}$. So this implies after base change a decomposition of the above $\mathcal{O}_F \otimes \mathcal{O}_{S_{K_N}}$-module $\text{Lie}(\mathcal{A})$.
But then according to this splitting and the above lemma, we get 
$\omega_{S_{K_N} \times F} = L_1 \oplus L_2,$
and therefore 
$\omega_{S_{K_N} \times F}^{\otimes 2} = L_1 \otimes L_2.$

For the $l$-adic Chern class of these line bundles we start with the Kummer
sequence of sheaves on $S_{K_N}\times\overline{\mathbb{Q}},$ that is
$0\rightarrow\boldsymbol{\mu}_{l^n} \rightarrow \mathbb{G}_{m} \rightarrow
\mathbb{G}_{m} \rightarrow 0.$ 
This gives rise to the Chern class map
$c_1 : \Pic(S_{K_N} \times \overline{\mathbb{Q}}) \rightarrow
H^2_{\acute{e}t}(S_{K_N} \times \overline{\mathbb{Q}},\mathbb{Q}_l(1)),$ 
where  $\Pic(S_{K_N} \times \overline{\mathbb{Q}}) = H^1_{\acute{e}t}(S_{K_N} \times \overline{\mathbb{Q}}, \mathbb{G}_{m,S_{K_N}}) = H^1_{\acute{e}t}(S_{K_N} \times \overline{\mathbb{Q}}, \mathcal{O}^*_{S_{K_N}})$.
\begin{Rmk}\label{RMKPOT}

To get rid of the level $N$-structure, we proceed as in the proof of
\cite{Ra}, Lemme 6.12, i.e. we look at the
$\left(K_N/K_{3N}\right)$-invariants in $\Pic(S_{K_{3N}} \times
\overline{\mathbb{Q}})$, resp. $H^2_{\acute{e}t}(S_{K_{3N}} \times
\overline{\mathbb{Q}},\mathbb{Q}_l(1))$. 
\end{Rmk}
Finally, we use our considerations in the complex case to get
\begin{Cor}\label{CORLCC}

Let $L_1 \otimes L_2$ be the line bundle of Hilbert modular forms of weight
$(2,2)$ on $S\times \overline{\mathbb{Q}}$. Then there is a unique line bundle $\widetilde{L_i}$, $i = 1,2$, in $\Pic(\widetilde{S}\times \overline{\mathbb{Q}} )$ with trivial Chern class on the boundary,
$\widetilde{L_i}|_{\widetilde{S}_{\infty} \times \overline{\mathbb{Q}}} \in
\text{Pic}^0(\widetilde{S}_{\infty} \times \overline{\mathbb{Q}}),$  such that
its restriction to the open part $S\times \overline{\mathbb{Q}}$ is $L_i$. In
particular, we have $c_1(L_i) \in H^2_{\acute{e}t,!}(S \times \overline{\mathbb{Q}},\mathbb{Q}_l(1)).$
\end{Cor}
\begin{proof}
We have discussed this for the complex surface $S(\mathbb{C})$.
By the comparison theorems this is also true for the algebraic classes
$c_1(L_i) \in H^2_{\acute{e}t,!}(S \times
\overline{\mathbb{Q}},\mathbb{Q}_l(1)).$ By Lemma \ref{LEMLIFT} we get this
extension in the complex context and in Lemma \ref{LEMCLB} we showed that the
restriction class is defined over $F$ and in particular  $\widetilde{L_i} \in  \Pic(\widetilde{S}\times F )$.
\end{proof}
\subsection{The Construction of the Motive}\label{CCHE} \label{Construction}
Our motive is an extension of a Tate motive by a Dirichlet motive and we get
such an extension as a piece of a long exact cohomology sequence. Our approach
is slightly different from Harder's general proposal in \cite{HaMM}, 1.4. There he considers a smooth closed subscheme $Y \subset X$ of a smooth projective scheme $X$ over $\Spec(\mathbb{Q})$. For such a pair one gets the long exact sequence
$$... \rightarrow H^i_c(U,\mathbb{Z}) \rightarrow H^i(U,\mathbb{Z}) \rightarrow H^i(\dot{\mathcal{N}}Y) \rightarrow ...,$$
where $U:=X - Y$ is the open part, and $\dot{\mathcal{N}}Y$ denotes the
punctured normal bundle of $Y$. Here in our situation the closed subscheme
$Y = \widetilde{S}_{\infty} \subset \widetilde{S} = X$ is not smooth. It is a
divisor with normal crossings. Let $j:S\hookrightarrow\widetilde{S}$ be the open embedding of $S/\mathbb{Q}$
into the toroidal compactification $\widetilde{S}/\mathbb{Q}$. We have an exact triangle
$j_!\mathbb{Q}\rightarrow \mathbf{R}j_*\mathbb{Q}\rightarrow
\mathbf{R}j_*\mathbb{Q}/j_!\mathbb{Q} \rightarrow j_!\mathbb{Q}[1]$ of complexes of sheaves on $\widetilde{S}$, i.e. a sequence in the derived category of  $\widetilde{S}$. This gives rise to the long exact sequence $...\rightarrow H^2(\widetilde{S}, j_!\mathbb{Q})\rightarrow
H^2(\widetilde{S},\mathbf{R}j_*\mathbb{Q})\rightarrow
H^2(\widetilde{S},\mathbf{R}j_*\mathbb{Q}/j_!\mathbb{Q})\rightarrow ... \: .$ 
By definition we have $H^i(\widetilde{S}, j_!\mathbb{Q})=H^i_c(S,\mathbb{Q})$ and
$H^i(\widetilde{S},\mathbf{R}j_*\mathbb{Q})=H^i(S,\mathbb{Q}).$ These are the motives as in \cite[1]{HaLNM}, \cite[1]{HaMM}. We discuss the
realisations in more detail below in Section \ref{MMwC}. For example there are
the \emph{Tate motives} $\mathbb{Q}(-n) = H^{2n}(\mathbb{P}^n,\mathbb{Q})$.
Now the first step of the construction is the 
\begin{Lem}\label{LEMBLC}

Let $...\rightarrow H^2(\widetilde{S}, j_!\mathbb{Q})\rightarrow H^2(\widetilde{S},\mathbf{R}j_*\mathbb{Q})\rightarrow H^2(\widetilde{S},\mathbf{R}j_*\mathbb{Q}/j_!\mathbb{Q})\rightarrow ...$
be the long exact sequence as above. 
Then one can break from this the exact sequence 
$$ 0 \rightarrow H^1(\widetilde{S},\mathbf{R}j_*\mathbb{Q}/j_!\mathbb{Q}) \rightarrow H^2_c(S,\mathbb{Q}) \rightarrow H^2(S,\mathbb{Q}) \rightarrow H^2(\widetilde{S},\mathbf{R}j_*\mathbb{Q}/j_!\mathbb{Q})\rightarrow 0.$$
\end{Lem}
\begin{proof}
By Poincar\'e duality  we just have to prove the vanishing of
 $H^1(S,\mathbb{Q}).$ This is stated for example in \cite{HLR}, Satz 1.9 (compare also
 \cite{SL2(O)}), but without any proof. In the Betti realisation the vanishing can be seen by the use of group cohomology for $\Gamma = \PSL_2(\mathcal{O}_F).$ One has
$H^1(S(\mathbb{C}),\mathbb{Q}) = H^1(\Gamma,\mathbb{Q}) = \Hom_{\text{Groups}}(\Gamma,\mathbb{Q}) \simeq \Hom_{\text{Ab}}(\Gamma/[\Gamma,\Gamma],\mathbb{Q}).$
The quotient $\Gamma^{\text{ab}} = \Gamma/[\Gamma,\Gamma]$ is finite (compare
\cite{Serre}, Th\'eor\`eme 3), and so $H^1(S(\mathbb{C}),\mathbb{Q})$ must be
trivial. By the comparison theorems this holds then also for the other cohomology theories.
\end{proof}
The complex $\mathbf{R}j_*\mathbb{Q}/j_!\mathbb{Q}$ lives on the boundary
$\widetilde{S}_{\infty}$ and is isomorphic to the complex
$i^*\mathbf{R}j_*\mathbb{Q},$ where $i: \widetilde{S}_{\infty} \hookrightarrow
\widetilde{S}$ denotes the closed embedding. We compute its cohomology in the
following lemma, which is just a special case of a theorem by Pink (see
e.g. \cite{HaLNM}, 2.2.10).
\begin{Prop}\label{Con}
Let $j:S\hookrightarrow\widetilde{S}$ be the open embedding of the
Hilbert modular surface $S$ into its toroidal compactification
$\widetilde{S}$. 
 Then we have 
$H^1(\widetilde{S},\mathbf{R}j_*\mathbb{Q}/j_!\mathbb{Q})\simeq \mathbb{Q}(0)$
and
$H^2(\widetilde{S},\mathbf{R}j_*\mathbb{Q}/j_!\mathbb{Q})\simeq \mathbb{Q}(-2).$
\end{Prop}
\begin{proof}
We denote by 
$\widetilde{S}_{\infty,i,i+1} := \widetilde{S}_{\infty,i}\cap \widetilde{S}_{\infty,i+1}$ 
the intersection of the two components $\widetilde{S}_{\infty,i}$ and
$\widetilde{S}_{\infty,i+1}$ with $i=0,..., n-1$. 
Now we distinguish the two cases. First, the point $P$ is smooth, i.e. $P \notin \widetilde{S}_{\infty,i,i+1}$ for all
$i$. Then we get for the fibre over $P$,
$$H^q(\mathbb{G}_m,\mathbb{Q})=
\begin{cases}
\mathbb{Q}(0), & q=0 \\
\mathbb{Q}(-1), & q=1 \\
0, & q=2. 
\end{cases}$$
If the point $P = P_{i,i+1}$ is not smooth, i.e. $P_{i,i+1} \in
\widetilde{S}_{\infty,i,i+1}$, we get for the fibre over $P_{i,i+1}$,
$$H^q(\mathbb{G}_m \times \mathbb{G}_m,\mathbb{Q})=
\begin{cases}
\mathbb{Q}(0), & q=0 \\
\mathbb{Q}(-1) \oplus \mathbb{Q}(-1) , & q=1 \\
\mathbb{Q}(-2), & q=2. 
\end{cases}$$
This leads us to 
$$R^qj_*\mathbb{Q}/j_!\mathbb{Q}=
\begin{cases} 
\mathbb{Q}(0)_{\widetilde{S}_{\infty}},& q=0 \\
\underset{i}{\bigoplus} \mathbb{Q}(-1)_{\widetilde{S}_{\infty,i}},& q=1 \\ 
\underset{i}{\bigoplus} \mathbb{Q}(-2)_{\widetilde{S}_{\infty,i,i+1}},& q=2,
\end{cases}$$
Now we put the cohomology classes $H^p(\widetilde{S},R^qj_*\mathbb{Q}/j_!\mathbb{Q})$ into a diagram, which describes the $E_2^{pq}$-term of the spectral sequence, i.e.
$E_2^{pq} = H^p(\widetilde{S},R^qj_*\mathbb{Q}/j_!\mathbb{Q})$ 
looks like
$$\xymatrix{{\underset{P_{i,i+1}}{\bigoplus}} \mathbb{Q}(-2) \ar[drr]^{d^{02}} & 0 &  \\
  {\underset{\widetilde{S}_{\infty,i}}{\bigoplus}}\mathbb{Q}(-1) \ar[drr]^{d^{01}} & 0 & {\underset{\widetilde{S}_{\infty,i}}{\bigoplus}} \mathbb{Q}(-2) \\
 {\mathbb{Q}}(0) & {\mathbb{Q}}(0) &  {\underset{\widetilde{S}_{\infty,i}}{\bigoplus}} \mathbb{Q}(-1)}$$
All the other entries are zero.
Since our polygon is a closed chain of $\mathbb{P}^1$'s, we conclude that the differential
$d^{02}: \underset{P_{i,i+1}}{\bigoplus} \mathbb{Q}(-2) \rightarrow \bigoplus_{\widetilde{S}_{\infty,i}} \mathbb{Q}(-2)$ must have rank $n-1$. To see this, we observe that $d^{02}$ is given by 
$$(P_{0,1}, ..., P_{n-1,0}) \mapsto (\widetilde{S}_{\infty,0} - \widetilde{S}_{\infty,1} , \widetilde{S}_{\infty,1} - \widetilde{S}_{\infty,2}, ...,  \widetilde{S}_{\infty,n-1} - \widetilde{S}_{\infty,0}),$$ 
i.e. can be represented by an $(n \times n)$-matrix like this
$\left(\begin{smallmatrix}1 & -1 &  & &  &  \\
  & 1 & -1 &  & &  \\ 
 & & \ddots & \ddots  & & \\  & & & & -1 & \\
  & &  & & 1 & -1\\
 -1 &  &  & & & 1
\end{smallmatrix} \right),$
which is of rank $n-1$. Therefore the kernel of $d^{02}$ is one copy of $\mathbb{Q}(-2)$. 
This gives us
$H^2(\widetilde{S},\mathbf{R}j_*\mathbb{Q}/j_!\mathbb{Q}) \simeq
\mathbb{Q}(-2).$ Furthermore, the kernel of $d^{21}$ is
$\text{Ker}(d^{21}) = \bigoplus_{\widetilde{S}_{\infty,i}}
\mathbb{Q}(-2)$. The quotient of this by the image of $d^{02}$ gives one copy
of $\mathbb{Q}(-2).$ And this is
$H^3(\widetilde{S},\mathbf{R}j_*\mathbb{Q}/j_!\mathbb{Q})$. The second differential map
$d^{01}:\bigoplus_{\widetilde{S}_{\infty,i}} \mathbb{Q}(-1)\rightarrow
\bigoplus_{\widetilde{S}_{\infty,i}} \mathbb{Q}(-1)$ is an isomorphism, hence we are left with
$H^1(\widetilde{S},\mathbf{R}j_*\mathbb{Q}/j_!\mathbb{Q}) \simeq \mathbb{Q}(0).$
\end{proof}
With the interior cohomology
$H^2_!(S,\mathbb{Q})=\text{Im}(H^2_c(S,\mathbb{Q})\rightarrow H^2(S,\mathbb{Q}))$
 we have an immediate consequence
\begin{Cor}\label{CC}
For a Hilbert modular surface $S$ we get a short exact sequence
$$0\rightarrow H^2_!(S,\mathbb{Q})\rightarrow H^2(S,\mathbb{Q})\rightarrow \mathbb{Q}(-2)\rightarrow 0.$$
\end{Cor}
Now we decompose the bottom of the extension $H^2_!(S,\mathbb{Q})$
even further. 
For this we twist the above sequence with $\mathbb{Q}(1)$ and get
$$0\rightarrow H^2_!(S,\mathbb{Q}(1)) \rightarrow H^2(S,\mathbb{Q}(1)) \rightarrow \mathbb{Q}(-1) \rightarrow 0.$$
By \cite{HLR}, 1.8, (see also \cite{vdG}, XI.2) we get, according to the semi-simple action of the Hecke algebra (see loc. cit. for the definition), a decomposition 
$H^2_!(S,\mathbb{Q}(1)) = H^2_{\cusp}(S,\mathbb{Q}(1)) \oplus H^2_{\CH}(S,\mathbb{Q}(1)),$
and furthermore, this Hecke action commutes with the action of Galois.
The first summand $H^2_{\cusp}(S,\mathbb{Q}(1))$ collects all the
contributions coming from the cuspidal representations of weight two. (In
\cite{HLR}, 1.8, this is denoted by $\text{Coh}_0$, the ``interesting'' part.)
The second one $H^2_{\CH}(S,\mathbb{Q}(1))$ consists of those coming from the
one dimensional representations, i.e. the Gr\"{o}\ss encharacters. (In
loc. cit. this is denoted by $\text{Coh}_e$, the ``trivial'' part.) Note that
in our case of $h^+=1$ and $K_f=K_0$ there is only a contribution by the
trivial character. Moreover, this summand $H^2_{\CH}(S,\mathbb{Q}(1))$ is
spanned by our two Chern classes $c_1(L_1)$ and $c_1(L_2)$ of Section
\ref{clb} (compare e.g. \cite{vdG}, XI.2).

Now we identify the Dirichlet character $\chi_D$ with the associated Galois character 
$\xymatrix{\Gal(\overline{\mathbb{Q}}/\mathbb{Q}) \ar[r]^>>>>>>{\Res} &
  \Gal(\mathbb{Q}(\zeta_D)/\mathbb{Q}) \ar[r]^>>>>>{\chi_D} &
  \Gal(F/\mathbb{Q}) \simeq \{\pm1\}},$ 
with the Galois group $\Gal(\mathbb{Q}(\zeta_D)/\mathbb{Q}) \simeq (\mathbb{Z}/D\mathbb{Z})^*$
of the cyclotomic field $\mathbb{Q}(\zeta_D)/\mathbb{Q}.$
Then define $\mathbb{Q}(0)\chi_D$ to be the Dirichlet motive for the quadratic character
$\chi_D$ in the sense of Deligne (\cite{De-Valuers}, 6). We explain this
latter notion in more detail below. And we get
\begin{Cor}\label{CORGAL}

We have an isomorphism
$$H^2_{\CH}(S,\mathbb{Q}(1)) \simeq \Res_{F/\mathbb{Q}}(\mathbb{Q}(0)) = \mathbb{Q}(0) \oplus \mathbb{Q}(0)\chi_D.$$
\end{Cor}
\begin{proof}
See e.g. \cite{HLR}, Proposition 2.10, or \cite{vdG}, XI.2 Proposition 2.7.
\end{proof}
If we use the above decomposition, 
 our sequence becomes 
$$0\rightarrow H^2_{\CH}(S,\mathbb{Q}(1)) \rightarrow H^2_{\CHE}(S,\mathbb{Q}(1)) \rightarrow \mathbb{Q}(-1)\rightarrow 0,$$
where $H^2_{\CHE}(S,\mathbb{Q}(1)):=H^2(S,\mathbb{Q}(1))/
H^2_{\cusp}(S,\mathbb{Q}(1)).$ So we identify the bottom $H^2_{\CH}(S,\mathbb{Q}(1))$ with $\mathbb{Q}(0) \oplus \mathbb{Q}(0)\chi_D$, and define
\begin{Def}[Kummer-Chern-Eisenstein Motive]

We call the extension 
$$0 \rightarrow \mathbb{Q}(0)\oplus \mathbb{Q}(0)\chi_D \rightarrow H^2_{\CHE}(S,\mathbb{Q}(1)) \rightarrow \mathbb{Q}(-1) \rightarrow 0$$
a Kummer-Chern-Eisenstein motive.
\end{Def}
Now we still have the action by the involution $\Theta \in
\text{Gal}(F/\mathbb{Q})$. This gives a further decomposition into $(\pm1)$-eigenspaces. Regarding this our sequence becomes after splitting the $(+1)$-eigenspace $\mathbb{Q}(0),$
$$0\rightarrow \mathbb{Q}(0)\chi_D \rightarrow H^2_{\CHE}(S,\mathbb{Q}(1))(-)\rightarrow \mathbb{Q}(-1) \rightarrow 0,$$
with $H^2_{\CHE}(S,\mathbb{Q}(1))(-):=H^2_{\CHE}(S,\mathbb{Q}(1))/\mathbb{Q}(0).$
So we get an element 
$$[H^2_{\CHE}(S,\mathbb{Q}(1))(-)] \in \Ext^1_{\mathcal{MM}_{\mathbb{Q}}}(\mathbb{Q}(-1),\mathbb{Q}(0)\chi_D),$$
which we call again a Kummer-Chern-Eisenstein motive, and which is now defined
over $\mathbb{Q}$. The name bases on the idea that it is an extension of $\mathbb{Q}(-1)$ by $\mathbb{Q}(0)\chi_D$ (so it should be Kummer), and that the extension $H^2_{\CHE}(S,\mathbb{Q}(1))$ is spanned by the Chern classes and the section of the restriction map. And the latter one is given by the Eisenstein section.

Now we have the various realisations of our motive. There are the $l$-adic realisations $[H^2_{\CHE,l}(S,\mathbb{Q}(1))(-)] \in
\Ext^1_{\mathcal{MGAL}}(\mathbb{Q}_l(-1),\mathbb{Q}_l(0)\chi_D),$ which are
mixed Galois modules (see e.g. \cite{HaLNM}, I). And we have the Hodge-de Rham realisation
$[H^2_{\CHE,\infty}(S,\mathbb{Q}(1))(-)] \in
\Ext^1_{\mathcal{MH}d\mathcal{R}_{\mathbb{Q}}}(\mathbb{Q}(-1),\mathbb{Q}(0)\chi_D),$
which is a mixed Hodge-de Rham structure (loc. cit.). 
\begin{Rmk}\label{RMKGEN}
We want to fix a generator of the summands of
$H^2_{\CH}(S,\mathbb{Q}(1))$. 
The first summand $\mathbb{Q}(0)$ is generated by $c_1(L_1 \otimes L_2) =
c_1(L_1) + c_1(L_2)$. 
For the second one $\mathbb{Q}(0)\chi_D$, we have the choice between $c_1(L_1
 \otimes L_2^{-1})$ and $c_1(L_1^{-1} \otimes L_2)$. These two differ
 by the action of $\Theta \in \Gal(F/\mathbb{Q})$. 
Here we choose $c_1(L_1 \otimes L_2^{-1})$ for a generator of $\mathbb{Q}(0)\chi_D$.
Note that we do not have a canonical choice of the generator, because both classes generate the submodule $\mathbb{Q}(0)\chi_D$. So we are left with this $(\pm 1)$-ambiguity. But our results respects this in the sense that if we flip the generator, we get the conjugate result.
\end{Rmk}
\subsection{The Dual Kummer-Chern-Eisenstein Motive} \label{dual motive}
In Chapter \ref{l-adic} we compute the $l$-adic realisations of $[H^2_{\CHE}(S,\mathbb{Q}(1))(-)]$. In order to do so, we look at the dual motive $[H^2_{\CHE}(S,\mathbb{Q}(1))(-)]^{\vee}$. 
\begin{Lem}\label{DM}

The Kummer-Chern-Eisenstein motive 
$$0 \rightarrow \mathbb{Q}(0) \chi_D \rightarrow H^2_{\CHE}(S,\mathbb{Q}(1))(-)\rightarrow \mathbb{Q}(-1) \rightarrow 0$$
becomes by dualising
$0\rightarrow \mathbb{Q}(1) \rightarrow
H^2_{\CHE}(S,\mathbb{Q}(1))(-)^{\vee}\rightarrow \mathbb{Q}(0) \chi_D
\rightarrow 0,$ where the bottom $\mathbb{Q}(1)$ is
$H^1(\widetilde{S},\mathbf{R}j_*\mathbb{Q}/j_!\mathbb{Q})\otimes
\mathbb{Q}(1)$ the motive as in Proposition \ref{Con}. In particular
$H^1(\widetilde{S},\mathbf{R}j_*\mathbb{Q}/j_!\mathbb{Q})$ is canonical dual
to $H^2(\widetilde{S},\mathbf{R}j_*\mathbb{Q}/j_!\mathbb{Q}).$
\end{Lem}
\begin{proof}

Recall that by Lemma \ref{LEMBLC} the sequence
$$ 0 \rightarrow H^1(\widetilde{S},\mathbf{R}j_*\mathbb{Q}/j_!\mathbb{Q}) \rightarrow H^2_c(S,\mathbb{Q}) \rightarrow H^2(S,\mathbb{Q}) \rightarrow H^2(\widetilde{S},\mathbf{R}j_*\mathbb{Q}/j_!\mathbb{Q})\rightarrow 0$$ 
is the starting point of the construction of the motive.
And the outer terms can be identified as  
$ 0 \rightarrow \mathbb{Q}(0) \rightarrow H^2_c(S,\mathbb{Q}) \rightarrow H^2(S,\mathbb{Q}) \rightarrow \mathbb{Q}(-2) \rightarrow 0,$
see Proposition \ref{Con}. To get the dual motive, we consider the Poincar\'e duality pairing
$H^2_c(S,\mathbb{Q}) \times H^2(S,\mathbb{Q}) \rightarrow \mathbb{Q}(-2).$ 
Hence $H^2_c(S,\mathbb{Q}(1))$ and $H^2(S,\mathbb{Q}(1))$ are dual to each other, i.e. the dual of the map 
$H^2_c(S,\mathbb{Q}(1)) \rightarrow H^2(S,\mathbb{Q}(1))$ 
is the map itself, and therefore the dual of the kernel is the cokernel and vice versa. Now this means that if we tensorise the above sequence with $\mathbb{Q}(1)$, its dual sequence is just the same sequence
$ 0 \rightarrow \mathbb{Q}(1) \rightarrow H^2_c(S,\mathbb{Q}(1)) \rightarrow
H^2(S,\mathbb{Q}(1)) \rightarrow \mathbb{Q}(-1) \rightarrow 0.$
So the dual sequence for 
$ 0  \rightarrow H^2_!(S,\mathbb{Q}(1)) \rightarrow H^2(S,\mathbb{Q}(1))
\rightarrow \mathbb{Q}(-1) \rightarrow 0$ 
is 
$ 0 \rightarrow \mathbb{Q}(1) \rightarrow H^2_c(S,\mathbb{Q}(1)) \rightarrow
\left[H^2_!(S,\mathbb{Q}(1))\right]^{\vee} \rightarrow 0.$ 
But the cokernel of $\mathbb{Q}(1) \rightarrow H^2_c(S,\mathbb{Q}(1))$ is 
$H^2_!(S,\mathbb{Q}(1))$  
the interior cohomology itself, this means $H^2_!(S,\mathbb{Q}(1))$ is
self-dual, that is
$ \left[H^2_!(S,\mathbb{Q}(1))\right]^{\vee} = H^2_!(S,\mathbb{Q}(1)).$ So take the Kummer-Chern-Eisenstein motive
$0 \rightarrow \mathbb{Q}(0)\chi_D \rightarrow H^2_{\CHE}(S,\mathbb{Q}(1))(-)
\rightarrow \mathbb{Q}(-1) \rightarrow 0.$ 
Then we get, according to the previous discussion, by dualising the extension 
$0 \rightarrow \mathbb{Q}(1) \rightarrow H^2_{\CHE}(S,\mathbb{Q}(1))(-)^{\vee}
\rightarrow \mathbb{Q}(0)\chi_D^{-1} \rightarrow 0.$ 
As $\chi_D^2 = \mathbf{1}$, we end up with 
$0 \rightarrow \mathbb{Q}(1) \rightarrow H^2_{\CHE}(S,\mathbb{Q}(1))(-)^{\vee} \rightarrow \mathbb{Q}(0)\chi_D \rightarrow 0$
in $\Ext^1_{\mathcal{MM}_{\mathbb{Q}}}(\mathbb{Q}(0)\chi_D,\mathbb{Q}(1))$.
\end{proof}
\begin{Rmk}\label{RMKDM}
The dual motive $H^2_{\CHE}(S,\mathbb{Q}(1))(-)^{\vee}$ sits in $H^2_c(S,\mathbb{Q}(1))$, because our motive
$[H^2_{\CHE}(S,\mathbb{Q}(1))(-)]$ is by definition a quotient of
$H^2(S,\mathbb{Q}(1))$ and $H^2(S,\mathbb{Q}(1))$ is the dual of
$H^2_c(S,\mathbb{Q}(1))$. We have chosen $c_1(L) = c_1(L_1\otimes L_2^{-1}) \in H^2_!(S,\mathbb{Q}(1))$ for the generator of $\mathbb{Q}(0)\chi_D$. 
In the dualising process we have to flip the generator and the dual motive of $\mathbb{Q}(0)\chi_D$ is therefore generated by
$c_1(L_1^{-1}\otimes L_2)$. Moreover, by Siegel's theorem (see
e.g. \cite[IV.1]{vdG} and Chapter \ref{HdR}) the cup product of $c_1(L_1\otimes L_2^{-1})$ with itself is
$-4\zeta_F(-1)$. And therefore to normalise the generator of the dual motive
we have to multiply it with $-\frac{1}{4\zeta_F(-1)}$.
\end{Rmk}
In the following we give an alternative construction of the dual extension. 
Start again with the open embedding $j$ of $S$ and the closed embedding $i$ of
$\widetilde{S}_{\infty}$ into the toroidal compactification $\widetilde{S}$. This gives a short exact sequence 
$0\rightarrow j_!j^*\mathbb{Q} \rightarrow \mathbb{Q} \rightarrow
i_*i^*\mathbb{Q} \rightarrow 0$ of sheaves on $\widetilde{S}$. 
By definition we have the identifications $H^i(\widetilde{S},j_!j^*\mathbb{Q}) = H^i_c(S,\mathbb{Q})$
and $H^i(\widetilde{S},i_*i^*\mathbb{Q}) =
H^i(\widetilde{S}_{\infty},\mathbb{Q}).$ 
Therefore, the sequence of sheaves induces a long exact sequence in the
cohomology. Since the compact surface $\widetilde{S}$ is simply connected
(cf. \cite{vdG}, IV.6 - more general the first cohomology with nontrivial
coefficients vanishes by  \cite{HLR}, Proposition 5.3), we get 
$$ 0 \rightarrow H^1(\widetilde{S}_{\infty},\mathbb{Q}) \rightarrow H^2_c(S,\mathbb{Q}) \stackrel{f_1}{\rightarrow} H^2(\widetilde{S},\mathbb{Q}) \rightarrow H^2(\widetilde{S}_{\infty},\mathbb{Q})\rightarrow 0.$$
We twist this with $\mathbb{Q}(1)$ and by the above remark we know that the dual motive $H^2_{\CHE}(S,\mathbb{Q}(1))(-)^{\vee}$ sits in $H^2_c(S,\mathbb{Q}(1))$. Moreover, we have 
\begin{Lem}\label{MOTKK}

Let 
$ 0 \rightarrow H^1(\widetilde{S}_{\infty},\mathbb{Q}(1)) \rightarrow
H^2_c(S,\mathbb{Q}(1)) \rightarrow \text{Im}(f_1) \rightarrow 0$ 
be the short exact sequence that is induced by the above sequence. Then the
dual of our Kummer-Chern-Eisenstein motive
$[H^2_{\CHE}(S,\mathbb{Q}(1))(-)]^{\vee}$ sits inside this extension, that is
$$\xymatrix{ 0 \ar[r] & H^1(\widetilde{S}_{\infty},\mathbb{Q}(1))  \ar@{=}[d] \ar[r] &  H^2_c(S,\mathbb{Q}(1)) \ar[r] &  {\text{Im}(f_1)} \ar[r] & 0 \\
 0 \ar[r] & { \mathbb{Q}(1)} \ar[r] & {H^2_{\CHE}(S,\mathbb{Q}(1))(-)^{\vee}} \ar[r] \ar@{^{(}->}[u] & {\mathbb{Q}(0) \chi_D} \ar[r] \ar@{^{(}->}[u] & 0.}$$
\end{Lem}
\begin{proof}
Consider the square 
$\xymatrix{H^2_c(S,\mathbb{Q}(1)) \ar[d]^{f_1} \ar@{=}[r] & H^2_c(S,\mathbb{Q}(1)) \ar[d]^{f_2} \\
 H^2(\widetilde{S},\mathbb{Q}(1)) \ar[r] &  H^2(S,\mathbb{Q}(1)).}$ 

By \cite{HLR}, Lemma 2.2 and Hilfssatz 2.3 (see also Proposition 5.3 of loc. cit.) we have that
$\text{Im}(f_2) = H^2_!(S,\mathbb{Q}(1)) \simeq \text{Im}(f_1)$
and moreover
$H^2_!(S,\mathbb{Q}(1)) \simeq \text{Im}\left(H^2(\widetilde{S},\mathbb{Q}(1)) \rightarrow H^2(S,\mathbb{Q}(1)) \right).$
So the two kernels 
$H^1(\widetilde{S}_{\infty},\mathbb{Q}(1)) = \text{Ker}(f_1)$ and 
$H^1(\widetilde{S},\mathbf{R}j_*\mathbb{Q}/j_!\mathbb{Q}) \otimes
\mathbb{Q}(1)  = \text{Ker}(f_2)$ have to be isomorphic, and we get the bottom
of our extension. Moreover, we have, by Lemma \ref{LEMLIFT}, the unique extension $c_1(\widetilde{L}) \in H^2(\widetilde{S},\mathbb{Q}(1))$  of $c_1(L) \in H^2_!(S,\mathbb{Q}(1))$ and $c_1(L)$ generates $\mathbb{Q}(0)\chi_D$. Since $\text{Im}(f_1) \simeq \text{Im}(f_2)$, we can consider $\mathbb{Q}(0)\chi_D$ inside $\text{Im}(f_1)$
\end{proof}
\subsection{Realisations of Mixed Motives}\label{MMwC}%
In this section we touch the general theory of motives a little
further. The main references are \cite{HaLNM}, 1, \cite{HaMM}, 1, and
\cite{De-Valuers}. (But consult additionally \cite{TK}, 1.)

According to the very general conjectures (see e.g. \cite{Ja}, 3) one expects for mixed Tate motives over a number field $k$ that
$\Ext^1_{\mathcal{MM}_{k}}(\mathbb{Q}(-1),\mathbb{Q}(0))= k^* \otimes \mathbb{Q}.$
But one expects even more. Each of such an extension should come from a so
called \emph{Kummer motive} $K \langle a \rangle$, $a \in k^*$. In other words
the Kummer motives should exhaust
$\Ext^1_{\mathcal{MM}_{k}}(\mathbb{Q}(-1),\mathbb{Q}(0))$. Recall briefly the construction of a Kummer motive $K \langle a
\rangle$ attached to $a \in k^*$. For details see e.g. \cite{TK}, 3.1. Start with $X = \mathbb{P}^1_k - \{\infty,0\} = \mathbb{G}_m$ and the divisor $D =\{1,a\}, \: a \neq 1$, and consider the relative cohomology sequence
$$0 \rightarrow H^0(\mathbb{G}_m,\mathbb{Q}) \rightarrow H^0(D,\mathbb{Q}) \rightarrow H^1(X,D,\mathbb{Q}) \rightarrow H^1(\mathbb{G}_m,\mathbb{Q}) \rightarrow 0,$$ 
which becomes 
$0 \rightarrow \mathbb{Q}(0) \rightarrow \mathbb{Q}(0) \oplus \mathbb{Q}(0) \rightarrow K \langle a \rangle \rightarrow \mathbb{Q}(-1) \rightarrow 0.$
And this gives us $K \langle a \rangle \in \Ext^1_{\mathcal{MM}_{k}}(\mathbb{Q}(-1),\mathbb{Q}(0)).$
The realisations are 
$(K \langle a \rangle_{\infty}, K \langle a \rangle_l) = (\prod_{\iota:k
  \hookrightarrow \mathbb{C}} \log |\iota(a)|, a).$ 
Now one observes (cf. Section \ref{OM}) that such a Kummer motive $K \langle a
\rangle$ is actually a  Kummer-$1$-motive $M_a$ in the sense of
\cite{De-HodgeIII}. Here we meet a slightly different situation, as our motive $[H^2_{\CHE}(S,\mathbb{Q}(1))(-)]$ is an extension of the Tate motive $\mathbb{Q}(-1)$ by the Dirichlet motive $\mathbb{Q}(0)\chi_D$. This means it is defined over $\mathbb{Q}$, but its extension class is actually in the real quadratic field $F$.
If we just look at the realisations of $[H^2_{\CHE}(S,\mathbb{Q}(1))(-)]$, the
situation becomes more transparent.
\subsubsection{The $l$-adic realisations}%
The $l$-adic realisations of our mixed motive
$[H^2_{\CHE}(S,\mathbb{Q}(1))(-)]$ is an extension of $\Gal(\overline{\mathbb{Q}}/\mathbb{Q})$-modules, i.e. an
element in $\Ext^1_{\mathcal{MGAL}}(\mathbb{Q}_l(-1),\mathbb{Q}_l(0)\chi_D).$ Let us describe this group. 
\begin{Lem}\label{LEMTMM}

Let $\Theta \in \Gal(F/\mathbb{Q})$ be the nontrivial element. Consider 
$\widehat{F^{*,(l)}} := \underset{n}{\varprojlim} \: {F}^*/\left(F^*\right)^{l^n}$ 
the $l$-adic completion of $F^*$ and define the subgroup of norm-one-elements
$\left(\widehat{F^{*,(l)}}\right)^{-\Theta} := \{f \in \widehat{F^{*,(l)}} | \Theta \cdot f = f^{-1} \}.$
Then there is a canonical isomorphism
$\Ext^1_{\mathcal{MGAL}}(\mathbb{Q}_l(-1),\mathbb{Q}_l(0)\chi_D) = \left(\widehat{F^{*,(l)}}\right)^{-\Theta} \otimes_{\mathbb{Z}_l} \mathbb{Q}_l.$
\end{Lem}
\begin{proof}

By definition we have 
$$\Ext^1_{\mathcal{MGAL}}(\mathbb{Q}_l(-1),\mathbb{Q}_l(0)\chi_D) = \underset{n}{\varprojlim} H^1(\Gal(\overline{\mathbb{Q}}/\mathbb{Q}),\boldsymbol{\mu}_{l^n} \otimes \chi_D) \otimes_{\mathbb{Z}_l} \mathbb{Q}_l.$$
Recall that $\boldsymbol{\mu}_{l^n} \otimes \chi_D$ denotes the Galois module,
given by the product of the cyclotomic character $\alpha$ and the quadratic
character $\chi_D$, i.e. the tensor product of Galois representations. Furthermore, we identify $\mathbb{Z}/l\mathbb{Z}$ with the roots of
unity $\boldsymbol{\mu}_l$. Now by the Hochschild-Serre spectral sequence one has the exact sequence
{\small$$0 \rightarrow H^1(\Gal(F/\mathbb{Q}),(\boldsymbol{\mu}_{l^n} \otimes \chi_D)^{\Gal(\overline{\mathbb{Q}}/F)}) \rightarrow {H^1({\Gal(\overline{\mathbb{Q}}/\mathbb{Q}}),\boldsymbol{\mu}_{l^n} \otimes \chi_D)}$$
$$ \rightarrow H^1(\Gal(\overline{\mathbb{Q}}/F),\boldsymbol{\mu}_{l^n} \otimes \chi_D)^{\Gal(F/\mathbb{Q})} \rightarrow H^2(\Gal(F/\mathbb{Q}),(\boldsymbol{\mu}_{l^n} \otimes \chi_D)^{\Gal(\overline{\mathbb{Q}}/F)}) \rightarrow. $$}
One firstly notes that $(\boldsymbol{\mu}_{l^n} \otimes \chi_D)^{\Gal(\overline{\mathbb{Q}}/F)}$ is trivial, i.e.
$$
H^1(\Gal(F/\mathbb{Q}),(\boldsymbol{\mu}_{l^n} \otimes \chi_D)^{\Gal(\overline{\mathbb{Q}}/F)})=0=  H^2(\Gal(F/\mathbb{Q}),(\boldsymbol{\mu}_{l^n} \otimes \chi_D)^{\Gal(\overline{\mathbb{Q}}/F)}).
$$ 
Hence the spectral sequence gives us
$${H^1({\Gal(\overline{\mathbb{Q}}/\mathbb{Q}}),\boldsymbol{\mu}_{l^n} \otimes \chi_D)} \simeq H^1(\Gal(\overline{\mathbb{Q}}/F),\boldsymbol{\mu}_{l^n} \otimes \chi_D)^{\Gal(F/\mathbb{Q})}.$$
Now $\Gal(\overline{\mathbb{Q}}/F)$ acts only on $\boldsymbol{\mu}_{l^n}$, and therefore we get a (twisted) Kummer isomorphism for $H^1(\Gal(\overline{\mathbb{Q}}/F),\boldsymbol{\mu}_{l^n} \otimes \chi_D)$.
This gives in the limit $\widehat{F^{*,(l)}} \otimes \chi_D$, where $\Gal(F/\mathbb{Q})$ acts on $\widehat{F^{*,(l)}}$ by conjugation, i.e. by $\Theta$, and on $\chi_D$ by the multiplication by $-1$. This means that, eventually, for the invariant part 
$$\left(\widehat{F^{*,(l)}} \otimes \chi_D \right)^{\Gal(F/\mathbb{Q})} \simeq \{f \in \widehat{F^{*,(l)}} | \Theta \cdot f = f^{-1} \}.$$
\end{proof}
So the $l$-adic realisations of the Kummer-Chern-Eisenstein motive is
determined by a scalar in $\widehat{F^{*,(l)}} \otimes_{\mathbb{Z}_l} \mathbb{Q}_l,$ which is
then independent of $l$. Moreover, the $l$-adic extension classes give rise to two dimensional Galois representations 
$\Gal(\overline{\mathbb{Q}}/\mathbb{Q}) \rightarrow \GL_2(\mathbb{Q}_l),$
with 
$\sigma \mapsto {\left(\begin{smallmatrix} \chi_D(\sigma) & * \\ 0 & \alpha^{-1}(\sigma) \end{smallmatrix}\right)},$
where $\alpha$ is the cyclotomic (Tate) character and where $*$ denotes the extension class. And by Kummer theory this element star $*$ is given by a Kummer field extension. If one suspects such an  $l$-adic representation to come from a Kummer motive $K \langle a \rangle$, with $a \in F^*$ of norm one, then the star $*$ is determined by the Kummer-one-cocycle $\frac{\sigma \cdot \left( \sqrt[l^{\infty}]{a}\right)}{\sqrt[l^{\infty}]{a}}$.   Hence $*= \tau_a(\sigma) \alpha^{-1}(\sigma)$, with 
$\frac{\sigma \cdot\left( \sqrt[l^{\infty}]{a}\right)}{\sqrt[l^{\infty}]{a}} = \zeta_{l^{\infty}}^{\tau_a(\sigma)}.$ 
Indeed we prove that our $[H^2_{\CHE,l}(S,\mathbb{Q}(1))(-)]$ induces such a
Galois representation (Theorem \ref{MTladic}).
\subsubsection{The Hodge-de Rham Realisation}
To describe
$\text{Ext}^1_{\mathcal{MH}d\mathcal{R}_{\mathbb{Q}}}(\mathbb{Q}(-1),\mathbb{Q}(0)\chi_D)$
we follow basically the exposition in \cite{TK}, 1.5 (see also \cite{HaMM},
1.5.2, and  \cite{HaLNM}, 4.3.2).

We are in the following situation: We have an element 
$[M]_{\infty} = [M]_{B-dR} \in \text{Ext}^1_{\mathcal{MH}d\mathcal{R}_{\mathbb{Q}}}(\mathbb{Q}(-1),\mathbb{Q}(0)\chi_D),$
i.e. two exact sequences in the commutative diagram
$$\xymatrix{ 0 \ar[r] & {\mathbb{Q}(0)\chi_{D,B}} \otimes \mathbb{C} \ar[d]^{\simeq}_{\cdot(\sqrt{D})^{-1}} \ar[r] & M_B \otimes \mathbb{C} \ar[d]^{\simeq}_{I_{\infty}} \ar[r] & {\mathbb{Q}(-1)_B} \otimes \mathbb{C} \ar[d]^{\simeq}_{\cdot (2 \pi i)^{-1}} \ar[r] \ar@/_1pc/[l]_{s_B} & 0 \\
0 \ar[r] & {\mathbb{Q}(0)\chi_{D,dR}} \otimes \mathbb{C} \ar[r] & M_{dR} \otimes \mathbb{C}\ar[r] & {\mathbb{Q}(-1)_{dR}} \otimes \mathbb{C} \ar[r] \ar@/^1pc/[l]^{s_{dR}} & 0, }$$ 
that are linked by the comparison isomorphism 
$I_{\infty} : M_B \otimes \mathbb{C} \simeq M_{dR} \otimes \mathbb{C}.$
The isomorphism on the right is given by multiplication with $ (2 \pi
i)^{-1}$, and on the left we have the multiplication with the inverse of the
Gau\ss-sum $G(\chi_D) = \sqrt{D}$ of $\chi_D$, see explanation below. So we get a representation like 
$\left(\begin{smallmatrix}  \sqrt{D} & * \\ 0 & 2 \pi i \end{smallmatrix}\right),$ 
where the star $*$ is our $[M]_{\infty}$. To get this, we observe that each of
the sequences splits in its own category, i.e. we have sections $s_B$ and
$s_{dR}$, and the extension class is given by the comparison of
$s_B$ and $s_{dR}$. Furthermore, we know for the Hodge filtration $F^{\bullet}M_{dR}$ of $M_{dR}$ that  
$F^{1}M_{dR} \simeq F^1 {\mathbb{Q}(-1)_{dR}} \simeq {\mathbb{Q}(-1)_{dR}},$
and this isomorphism gives the section  $s_{dR}$. Now let us describe the rule to get $[M]_{\infty} \in
\text{Ext}^1_{\mathcal{MH}d\mathcal{R}_{\mathbb{Q}}}(\mathbb{Q}(-1),\mathbb{Q}(0)\chi_D)$.
We start with a generator $\mathbf{1}_B$ of ${\mathbb{Q}(-1)_B} \otimes
\mathbb{C}$. This goes via the right isomorphism to $ (2 \pi i)^{-1} \cdot
\mathbf{1}_{dR} \in  \mathbb{Q}(-1)_{dR} \otimes \mathbb{C}$, by  $s_{dR}$ we
land in $M_{dR} \otimes \mathbb{C}$ and by $I_{\infty}^{-1}$ in $M_B \otimes
\mathbb{C}$. On the other hand we can map  $\mathbf{1}_B \in
{\mathbb{Q}(-1)_B} \otimes \mathbb{C}$ by the section  $s_B$ directly to $M_B
\otimes \mathbb{C}$ such that the image is in the $(-1)$-eigenspace of
$F_{\infty}$. Then the difference $s_B(\mathbf{1}_B) - I_{\infty}^{-1}( (2 \pi i)^{-1} \cdot
\mathbf{1}_{dR})$ is in the kernel $\mathbb{Q}(0)\chi_{D,B} \otimes
\mathbb{C}$. More precisely, it is in $i\mathbb{R}$. If we multiply this with
the inverse of the Gau\ss-sum $G(\chi_D)$ of $\chi_D$ (here it is
$(\sqrt{D})^{-1}$), we get exactly our class $[M]_{\infty}$. We summarise the description of the Hodge-de Rham extension classes in
\begin{Lem}\label{LEMHSIP} \label{RMKBASE}

There is an isomorphism 
$\text{Ext}^1_{\mathcal{MH}d\mathcal{R}_{\mathbb{Q}}}(\mathbb{Q}(-1),\mathbb{Q}(0)\chi_D) \simeq i\mathbb{R}.$
We choose here $\frac{1}{2 \pi i}$ as a basis for $i\mathbb{R}$, i.e. we get 
$\text{Ext}^1_{\mathcal{MH}d\mathcal{R}_{\mathbb{Q}}}(\mathbb{Q}(-1),\mathbb{Q}(0)\chi_D) \simeq \mathbb{R}.$
\end{Lem}
\begin{proof}

For details see e.g. \cite{TK}, 1.5.
\end{proof}
Note that the image of the realisation functor consists of those elements in
$\mathbb{R}$ such that $\Theta$ acts by $-1$. Again in the case of a Kummer motive $K \langle a \rangle$ one knows how the
realisation class  $K \langle a \rangle_{\infty}$ looks like (compare
\cite{TK}, 3.1), one gets 
$\left(\begin{smallmatrix}  \sqrt{D} & \log a \\ 0 & 2 \pi i \end{smallmatrix}\right).$ So in view of the above conjecture, one should have the following identification
$\Ext^1_{\mathcal{MM}_{\mathbb{Q}}}(\mathbb{Q}(-1),\mathbb{Q}(0)\chi_D) \simeq
(F^*)^{-\Theta} \otimes \mathbb{Q},$ 
i.e. if we have got a pair $(M_{\infty}, M_l)$ such that there exists an $a \in F^*$ of norm one with
$(M_{\infty}, M_l) = (\log a,a),$
then this should come from a Kummer motive $K \langle a \rangle \in \Ext^1_{\mathcal{MM}_{\mathbb{Q}}}(\mathbb{Q}(-1),\mathbb{Q}(0)\chi_D)$.

In the end of this section I discuss briefly the category of
motives, where the constructed Chern-Eisenstein motives are living in. It is
the category that is generated by Dirichlet motives and their extensions -
cf. \cite{De-Valuers}, 6. One should think of a Dirichlet motive in our special case of a real quadratic
number field $F/\mathbb{Q}$ and its character $\chi_D$ as follows: our $F/\mathbb{Q}$ is a subfield of $\mathbb{Q}(\zeta_D)/\mathbb{Q}$, where $D$ is the discriminant of $F$. Hence we get 
$\Spec(\mathbb{Q}(\zeta_D)) \rightarrow \Spec(F).$
The abelian Galois module $\Spec(\mathbb{Q}(\zeta_D))$ decomposes by the characters, and we define $\mathbb{Q}(0)\chi_D$ as the direct summand, corresponding to the projector 
$\frac{1}{|(\mathbb{Z}/D\mathbb{Z})^*|} \sum_{\sigma \in
  \Gal(\mathbb{Q}(\zeta_D)/\mathbb{Q})} \chi_D(\sigma) \cdot \sigma,$
 where we use that $\chi_D^{-1} = \chi_D$. (The general definition of
 $\mathbb{Q}(0)\chi$ for an arbitrary Dirichlet character $\chi :
 (\mathbb{Z}/D\mathbb{Z})^* \rightarrow E^*$ with values in a field $E$ is
 given by the above projector with the inverse character $\chi^{-1}$.) If we
 restrict this to $\Spec(F)$, we get  $\mathbb{Q}(0)\chi_D$ as the direct
 summand, corresponding to the projector $\frac{1}{2}(\text{id}-\Theta)$,
 where $\Theta$ is the nontrivial element in $\Gal(F/\mathbb{Q})$. So we have
 $\Spec(F)$ as a two dimensional motive over $\Spec(\mathbb{Q})$, i.e. $H^0(\Spec(F)) = \mathbb{Q}(0) \oplus \mathbb{Q}(0)\chi_D.$

Let us again look at the realisations of $\mathbb{Q}(0)\chi_D$. There we have
(by loc. cit.)
the $l$-adic realisations that are one dimensional Galois modules
$\mathbb{Q}_l(0)\chi_D$, i.e. $\Gal(\overline{\mathbb{Q}}/\mathbb{Q})$ acts by
$\chi_D$. For the Hodge-de Rham realisation we have that the Betti realisation
$\mathbb{Q}(0)\chi_{D,B}$ is just $\mathbb{Q}(0)$. The Hodge structure of
$\mathbb{Q}(0) \otimes \mathbb{C}$ is pure of type $(0,0)$, and the involution
$F_{\infty}$ given by the complex conjugation. Furthermore, over $\mathbb{C}$
it becomes $\mathbb{Q}(0)\chi_{D,B} \otimes \mathbb{C}$ and isomorphic to the
de Rham realisation $\mathbb{Q}(0)\chi_{D,dR} \otimes \mathbb{C}$. As an
appropriate basis of $\mathbb{Q}(0)\chi_{D,dR}$ we choose the Gau\ss-sum
$G(\chi_D)$ of $\chi_D$ - compare \cite{De-Valuers}, 6.4. In our case of the
quadratic character (with $D \equiv 1 (\text{mod} \: 4)$) this is $G(\chi_D) =
\sqrt{D}$. But this is a very special case of a Dirichlet motive, in particular its values are yet in $\mathbb{Q}$.

\section{The $l$-adic Realisations of our Motive}\label{l-adic}
We start with the calculation of the $l$-adic realisations of
our Kummer-Chern-Eisenstein motive. We determine the action of
$\Gal(\overline{\mathbb{Q}}/\mathbb{Q})$ on
$[H^2_{\CHE,l}(S,\mathbb{Q}(1))(-)]$.  

For each $l$ we get a two dimensional extension
$$0\rightarrow  \mathbb{Q}_l(0)\chi_D \rightarrow H^2_{\CHE,l}(S,\mathbb{Q}(1))(-) \rightarrow \mathbb{Q}_l(-1) \rightarrow 0,$$
 and we know that $\Gal(\overline{\mathbb{Q}}/\mathbb{Q})$ acts by the
quadratic $\chi_D$ on the left and by the inverse of the Tate character
 $\alpha$ on the right term. Hence there is a two dimensional Galois representation $\Gal(\overline{\mathbb{Q}}/\mathbb{Q}) \rightarrow \GL_2(\mathbb{Q}_l),$
with 
$\sigma \mapsto {\left(\begin{smallmatrix}  \chi_D(\sigma) & * \\  0 & \alpha^{-1}(\sigma) \end{smallmatrix}\right)}$
and the star $*$ represents the extension class in
$\Ext^1_{\mathcal{MGAL}}(\mathbb{Q}_l(-1),\mathbb{Q}_l(0)\chi_D).$ Now in
accordance of Lemma \ref{LEMTMM} this group is $\left(\widehat{F^{*,(l)}}\right)^{-\Theta} \otimes_{\mathbb{Z}_l} \mathbb{Q}_l.$
In the following sections I compute this scalar, which belongs to the $l$-adic realisations of the Kummer-Chern-Eisenstein motive. We know by Kummer theory that this gives rise to Kummer field extensions of $F$, and the associated Galois representations of our $l$-adic realisations above come exactly from this field extensions.
\subsection{Detection of the $l$-adic Extension Classes} \label{DoC}%
The first problem is to find an appropriate recipe that describes such an extension class. In our situation this can be done in the following manner.

If we look at the dual motive (compare Lemma \ref{DM}) we get extensions
$$0\rightarrow \mathbb{Q}_l(1) \rightarrow H^2_{\CHE,l}(S,\mathbb{Q}(1))(-)^{\vee} \rightarrow \mathbb{Q}_l(0)\chi_D \rightarrow 0$$
in $\Ext^1_{\mathcal{MGAL}}(\mathbb{Q}_l(0)\chi_D,\mathbb{Q}_l(1)).$ Hence the Galois representation for the dual module looks like
$ \sigma \mapsto \left(\begin{smallmatrix} \alpha(\sigma) & * \\ 0 &  \chi_D(\sigma) \end{smallmatrix}\right).$
Recall by Section \ref{Construction} and \ref{dual motive} that the top of the extension $\mathbb{Q}_l(0)\chi_D$ is generated by the first Chern class 
$c_1(L) = c_1(L_1^{-1} \otimes L_2) \in H^2_{\acute{e}
  t,!}(S\times\overline{\mathbb{Q}},\mathbb{Q}_l(1)) \subset H^2_{\acute{e}
  t}(S\times\overline{\mathbb{Q}},\mathbb{Q}_l(1)),$ and note that we write again $L$ (instead of $L^{-1}$)
to keep the notation easy and we normalise as described in Remark
\ref{RMKDM}. Furthermore, the middle
$H^2_{\CHE,l}(S,\mathbb{Q}(1))(-)^{\vee}$ sits in the cohomology with compact
support $H^2_{\acute{e} t,c}(S\times\overline{\mathbb{Q}},\mathbb{Q}_l(1))$ -
see loc. cit. And the bottom $\mathbb{Q}_l(1)$ comes from the cohomology $H^1_{\acute{e} t}(\widetilde{S}\times\overline{\mathbb{Q}},\mathbf{R}j_*\mathbb{Q}_l/j_!\mathbb{Q}_l)$ twisted by $\mathbb{Q}_l(1)$.

Now along the general procedure to get such a group extension class, we must
lift $c_1(L) \in H^2_{\acute{e}
  t,!}(S\times\overline{\mathbb{Q}},\mathbb{Q}_l(1))$ to $\widetilde{c_1(L)}
\in H^2_{\acute{e} t,c}(S\times\overline{\mathbb{Q}},\mathbb{Q}_l(1))$ and
find the action of $\Gal(\overline{\mathbb{Q}}/\mathbb{Q})$ on this lifting
$\widetilde{c_1(L)}$. The extension class is then given by the cocycle
$\sigma \left(\widetilde{c_1(L)}\right) - \widetilde{\sigma c_1(L)}, \: \sigma \in \Gal(\overline{\mathbb{Q}}/\mathbb{Q}).$
This is an element in $H^2_{\acute{e}
  t,c}(S\times\overline{\mathbb{Q}},\mathbb{Q}_l(1))$ that goes to zero in
$H^2_{\acute{e} t}(S\times\overline{\mathbb{Q}},\mathbb{Q}_l(1)),$ i.e. it is
in the kernel, which is $\mathbb{Q}_l(1)$. 
We note that the above cocycle $\sigma \left(\widetilde{c_1(L)}\right) -
\widetilde{\sigma c_1(L)}$ is in $\Ext^1_{\mathcal{MGAL}}(\mathbb{Q}_l(0)\chi_D,\mathbb{Q}_l(1)).$ 
This is just because $\Theta \in \Gal(F/\mathbb{Q})$ flips the two factors of
$L = L_1^{-1} \otimes L_2$, i.e. $\Theta \cdot c_1(L) = -  c_1(L)$.
We observe that we have the first Chern class of a line bundle, which is defined over $F$, i.e. it is in the Galois invariant part $H^2_{\acute{e} t}(S\times\overline{\mathbb{Q}},\mathbb{Q}_l(1))^{\Gal(\overline{\mathbb{Q}}/F)}$. This means that the Galois group $\Gal(\overline{\mathbb{Q}}/F)$ acts trivially on $c_1(L)$, and in particular, 
$\sigma c_1(L) = c_1(L)$, resp. $\widetilde{\sigma c_1(L)} = \widetilde{c_1(L)}.$ Recall that (by Section \ref{clb}) the above $c_1(L_i)$'s in $H^2_{\acute{e} t,!}(S\times\overline{\mathbb{Q}},\mathbb{Q}_l(1))$ come uniquely from the first Chern classes of the line bundles $\widetilde{L_i}$ of modular forms on the toroidal compactified surface $\widetilde{S}\times\overline{\mathbb{Q}}$, i.e.
$c_1(L_i) = c_1(\widetilde{L_i})|_{S \times \overline{\mathbb{Q}}}.$ So in
order to lift our class $c_1(L) \in H^2_{\acute{e}
  t,!}(S\times\overline{\mathbb{Q}},\mathbb{Q}_l(1))$ to $H^2_{\acute{e}
  t,c}(S\times\overline{\mathbb{Q}},\mathbb{Q}_l(1))$, we can lift the class
$c_1(\widetilde{L})$ in $H^2_{\acute{e}
  t}(\widetilde{S}\times\overline{\mathbb{Q}},\mathbb{Q}_l(1))$ to
$H^2_{\acute{e} t,c}(S\times\overline{\mathbb{Q}},\mathbb{Q}_l(1)).$ The following lemma shows that those two liftings are equal. 
\begin{Lem}\label{LEMLL}

Start with $c_1(L) \in H^2_{\acute{e}
  t,!}(S\times\overline{\mathbb{Q}},\mathbb{Q}_l(1))$ and lift this to $\widetilde{c_1(L)} \in H^2_{\acute{e}
  t,c}(S\times\overline{\mathbb{Q}},\mathbb{Q}_l(1))$. Let
$c_1(\widetilde{L}) \in H^2_{\acute{e}
  t}(\widetilde{S}\times\overline{\mathbb{Q}},\mathbb{Q}_l(1))$ be the Chern
class of the extended line bundle $\widetilde{L}.$ Then $\widetilde{c_1(L)}$ maps to $c_1(\widetilde{L})$ via $H^2_{\acute{e} t,c}(S\times\overline{\mathbb{Q}},\mathbb{Q}_l(1)) \rightarrow H^2_{\acute{e} t}(\widetilde{S}\times\overline{\mathbb{Q}},\mathbb{Q}_l(1))$, i.e.  $c_1(\widetilde{L})$ lifts to $\widetilde{c_1(L)}$.
\end{Lem}
\begin{proof}
Again according to \cite{HLR}, Lemma 2.2, we know that $H^2_{\acute{e}
  t,!}(S\times\overline{\mathbb{Q}},\mathbb{Q}_l(1))$ is isomorphic to the
  image of the restriction map $H^2_{\acute{e}
  t}(\widetilde{S}\times\overline{\mathbb{Q}},\mathbb{Q}_l(1)) \rightarrow
  H^2_{\acute{e} t}(S\times\overline{\mathbb{Q}},\mathbb{Q}_l(1)).$ 
And moreover, we have that $c_1(\widetilde{L}) \in H^2_{\acute{e}
  t}(\widetilde{S}\times\overline{\mathbb{Q}},\mathbb{Q}_l(1))$ is the unique
lift of $c_1(L) \in H^2_{\acute{e}
  t}(S\times\overline{\mathbb{Q}},\mathbb{Q}_l(1))$ - see Corollary \ref{CORLCC}.
\end{proof} 
Our next goal is the construction of a diagram, which contains all cohomology
 groups that play along. For this we turn over to finite coefficient $\boldsymbol{\mu}_{l^n}$, i.e. we
 consider the map $c_1^{(l^n)}$ that is given by
$\xymatrix{ H^1_{\acute{e} t}(\widetilde{S}\times\overline{\mathbb{Q}},\mathbb{G}_m) \ar[r]^{c_1} \ar[dr]_{c_1^{(l^n)}} & H^2_{\acute{e} t}(\widetilde{S}\times\overline{\mathbb{Q}},\mathbb{Z}_{l}(1)) \ar[d]^{c_1 \text{mod} \: l^n} \\
 &  H^2_{\acute{e}
   t}(\widetilde{S}\times\overline{\mathbb{Q}},\boldsymbol{\mu}_{l^n})}$ 
and we get the class 
$c_1^{(l^n)}(\widetilde{L}) \in H^2_{\acute{e} t}(\widetilde{S}\times\overline{\mathbb{Q}},\boldsymbol{\mu}_{l^n})$,
and furthermore we observe that
$c_1(\widetilde{L}) = \underset{n}{\varprojlim} \: c_1^{(l^n)}(\widetilde{L})
\in  H^2_{\acute{e}
  t}(\widetilde{S}\times\overline{\mathbb{Q}},\mathbb{Z}_{l}(1)).$ 
We know (Lemma \ref{DM}) that the bottom $\mathbb{Q}_l(1)$ of our extension comes from $H^1_{\acute{e} t}(\widetilde{S}\times\overline{\mathbb{Q}},\mathbf{R}j_*\mathbb{Q}_l/j_!\mathbb{Q}_l) \otimes \mathbb{Q}_l(1)$, 
and the $l$-adic version of our Lemma \ref{MOTKK} gives the isomorphism 
 $\underset{n}{\varprojlim}H^1_{\acute{e}
   t}(\widetilde{S}_{\infty}\times\overline{\mathbb{Q}},\boldsymbol{\mu}_{l^n}) \otimes \mathbb{Q}_l = H^1_{\acute{e} t}(\widetilde{S}\times\overline{\mathbb{Q}},\mathbf{R}j_*\mathbb{Q}_l/j_!\mathbb{Q}_l) \otimes \mathbb{Q}_l(1).$ 
To construct the diagram, we start with the open embedding $j$ of $S \times \overline{\mathbb{Q}}$ and the
closed embedding $i$ of $\widetilde{S}_{\infty} \times \overline{\mathbb{Q}}$
into the toroidal compactification $\widetilde{S} \times
\overline{\mathbb{Q}}$. 
We have the following diagram of sheaves on $\widetilde{S} \times \overline{\mathbb{Q}}$ with exact rows and columns.
$$\xymatrix{ & 0 \ar[d] & 0 \ar[d] & 0 \ar[d] & \\
0 \ar[r] & {j_!j^*\boldsymbol{\mu}_{l^n,\widetilde{S}}} \ar[d] \ar[r]\ar[dr]_{\alpha'} & {j_!j^*\mathbb{G}_{m,\widetilde{S}}} \ar[d] \ar[r] & {j_!j^*\mathbb{G}_{m,\widetilde{S}}} \ar[d] \ar[r] & 0 \\
0 \ar[r] & {\boldsymbol{\mu}_{l^n,\widetilde{S}}}  \ar[d] \ar[r]^{\alpha_1} & {\mathbb{G}_{m,\widetilde{S}}} \ar[d] \ar[r]\ar[dr]_{\beta} & {\mathbb{G}_{m,\widetilde{S}}} \ar[d] \ar[r] & 0 \\
0 \ar[r] &  i_*i^*\boldsymbol{\mu}_{l^n,\widetilde{S}} \ar[d] \ar[r] & {i_*i^*\mathbb{G}_{m,\widetilde{S}}} \ar[d] \ar[r] & {i_*i^*\mathbb{G}_{m,\widetilde{S}}} \ar[d] \ar[r]& 0\\ 
& 0 & 0 & 0 & }$$
where in each row sits a Kummer sequence. 
Now take a look at the cohomology of these sheaves, i.e. we have
\begin{Lem} \label{BD}

The above diagram of sheaves on $\widetilde{S} \times \overline{\mathbb{Q}}$ gives rise to the following diagram of cohomology groups with Galois equivariant maps and exact rows and columns. 
{\scriptsize $$\xymatrix@C=0.8em{ && & 0 \ar[d]\ar[r] & H^1_{\acute{e} t}(\widetilde{S}\times\overline{\mathbb{Q}},\mathbb{G}_m) \ar[d] \\
& & & H^1_{\acute{e} t}(\widetilde{S}_{\infty}\times\overline{\mathbb{Q}},\boldsymbol{\mu}_{l^n}) \ar[d] \ar[r] & H^1_{\acute{e} t}(\widetilde{S}_{\infty}\times\overline{\mathbb{Q}},i^*\mathbb{G}_{m,\widetilde{S}}) \ar[d] \\  
& H^1_{\acute{e} t,c}(S\times\overline{\mathbb{Q}},\mathbb{G}_m) \ar[d] \ar[r] & H^1_{\acute{e} t,c}(S\times\overline{\mathbb{Q}},\mathbb{G}_m) \ar[d] \ar[r] & H^2_{\acute{e} t,c}(S\times\overline{\mathbb{Q}},\boldsymbol{\mu}_{l^n}) \ar[d] \ar[r] & H^2_{\acute{e} t,c}(S\times\overline{\mathbb{Q}},\mathbb{G}_m) \ar[d] \\
0 \ar[r] \ar[d] & H^1_{\acute{e} t}(\widetilde{S}\times\overline{\mathbb{Q}},\mathbb{G}_m) \ar[d] \ar[r] & H^1_{\acute{e} t}(\widetilde{S}\times\overline{\mathbb{Q}},\mathbb{G}_m) \ar[d] \ar[r] & H^2_{\acute{e} t}(\widetilde{S}\times\overline{\mathbb{Q}},\boldsymbol{\mu}_{l^n})\ar[d] \ar[r] & H^2_{\acute{e} t}(\widetilde{S}\times\overline{\mathbb{Q}},\mathbb{G}_m) \\
 H^1_{\acute{e} t}(\widetilde{S}_{\infty}\times\overline{\mathbb{Q}},\boldsymbol{\mu}_{l^n}) \ar[d] \ar[r] & H^1_{\acute{e} t}(\widetilde{S}_{\infty}\times\overline{\mathbb{Q}},i^*\mathbb{G}_{m,\widetilde{S}}) \ar[r] \ar[d] & H^1_{\acute{e} t}(\widetilde{S}_{\infty}\times\overline{\mathbb{Q}},i^*\mathbb{G}_{m,\widetilde{S}}) \ar[r]  & H^2_{\acute{e} t}(\widetilde{S}_{\infty}\times\overline{\mathbb{Q}},\boldsymbol{\mu}_{l^n}) \ar[d] & \\ H^2_{\acute{e} t,c}(S\times\overline{\mathbb{Q}},\boldsymbol{\mu}_{l^n}) \ar[r] & H^2_{\acute{e} t,c}(S\times\overline{\mathbb{Q}},\mathbb{G}_m)  & & 0 & }$$}
\end{Lem}
\begin{proof}
By definition we have  
$H^i_{\acute{e}
  t}(\widetilde{S}\times\overline{\mathbb{Q}},j_!j^*\mathbb{G}_{m,\widetilde{S}}) = H^i_{\acute{e} t,c}(S\times\overline{\mathbb{Q}},\mathbb{G}_m)$ 
and 
$H^i_{\acute{e}
  t}(\widetilde{S}\times\overline{\mathbb{Q}},j_!j^*\boldsymbol{\mu}_{l^n,\widetilde{S}}) = H^i_{\acute{e} t,c}(S\times\overline{\mathbb{Q}},\boldsymbol{\mu}_{l^n}).$ 
For $\boldsymbol{\mu}_{l^n}$ one also knows 
$H^i_{\acute{e}
  t}(\widetilde{S}\times\overline{\mathbb{Q}},i_*i^*\boldsymbol{\mu}_{l^n,\widetilde{S}}) = H^i_{\acute{e} t}(\widetilde{S}_{\infty}\times\overline{\mathbb{Q}},\boldsymbol{\mu}_{l^n}).$ 
But for $\mathbb{G}_m$ we only know 
$H^i_{\acute{e} t}(\widetilde{S}\times\overline{\mathbb{Q}},i_*i^*\mathbb{G}_{m,\widetilde{S}}) \simeq H^i_{\acute{e} t}(\widetilde{S}_{\infty} \times \overline{\mathbb{Q}},i^*\mathbb{G}_{m,\widetilde{S}}).$  The compact surface $\widetilde{S}$ is simply connected
(cf. \cite{vdG}, IV.6), i.e. $H^1_{\acute{e}
  t}(\widetilde{S}\times\overline{\mathbb{Q}},\boldsymbol{\mu}_{l^n})$
vanishes (cf. also \cite{HLR}, Proposition 5.3).
\end{proof} 
The $l$-adic version of Lemma \ref{MOTKK} tells us that we can find our dual extension in the above big diagram.
Let 
$$0 \rightarrow  \underset{n}{\varprojlim}H^1_{\acute{e} t}(\widetilde{S}_{\infty} \times\overline{\mathbb{Q}},\boldsymbol{\mu}_{l^n}) \otimes \mathbb{Q}_l  \underset{f_1}{\rightarrow}   \underset{n}{\varprojlim}H^2_{\acute{e} t,c}(S\times\overline{\mathbb{Q}},\boldsymbol{\mu}_{l^n}) \otimes \mathbb{Q}_l  \rightarrow \text{Im}(f_1) \rightarrow 0$$ 
be the short exact sequence, coming from the above diagram in Lemma \ref{BD}.
Then the $l$-adic realisation $[H^2_{\CHE,l}(S,\mathbb{Q}(1))(-)]^{\vee}$ sits inside this extension.
\subsection{Galois Action on the Liftings}\label{liftings}%
The idea to describe the Galois action is based on the above big diagram in Lemma \ref{BD}. It relates the Chern class $c_1^{(l^n)}(\widetilde{L}) \in H^2_{\acute{e} t}(\widetilde{S}\times\overline{\mathbb{Q}},\boldsymbol{\mu}_{l^n})$  of the line bundle $\widetilde{L} \in H^1_{\acute{e} t}(\widetilde{S}\times\overline{\mathbb{Q}},\mathbb{G}_m)$ and the restriction of $\widetilde{L}$ to the boundary $\widetilde{S}_{\infty}\times\overline{\mathbb{Q}}$.  

To get this restriction $\widetilde{L}|_{\widetilde{S}_{\infty} \times \overline{\mathbb{Q}}}$, we compose the above ``restriction map'' in the diagram 
$R: H^1_{\acute{e} t}(\widetilde{S}\times\overline{\mathbb{Q}},\mathbb{G}_m) \rightarrow H^1_{\acute{e} t}(\widetilde{S}_{\infty}\times\overline{\mathbb{Q}},i^*\mathbb{G}_{m,\widetilde{S}})$ 
with the map 
$H^1_{\acute{e} t}(i^*):H^1_{\acute{e} t}(\widetilde{S}_{\infty}\times\overline{\mathbb{Q}},i^*\mathbb{G}_{m,\widetilde{S}}) \rightarrow H^1_{\acute{e} t}(\widetilde{S}_{\infty}\times\overline{\mathbb{Q}},\mathbb{G}_{m}),$
which is induced by the morphism $i^*\mathbb{G}_{m,\widetilde{S}} \rightarrow
\mathbb{G}_{m,\widetilde{S}_{\infty}}$ of sheaves on
$\widetilde{S}_{\infty}\times\overline{\mathbb{Q}}$. Now we add to the bottom
part of the diagram in Lemma \ref{BD} the map $H^1_{\acute{e} t}(i^*)$, that
is 
$$\xymatrix{ & H^1_{\acute{e} t}(\widetilde{S}\times\overline{\mathbb{Q}},\mathbb{G}_m) \ar[d]^{R} \ar[r] & H^1_{\acute{e} t}(\widetilde{S}\times\overline{\mathbb{Q}},\mathbb{G}_m) \ar[d]^{R} \\
H^1_{\acute{e} t}(\widetilde{S}_{\infty}\times\overline{\mathbb{Q}},\boldsymbol{\mu}_{l^n}) \ar[r] \ar@{=}[d] & H^1_{\acute{e} t}(\widetilde{S}_{\infty}\times\overline{\mathbb{Q}},i^*\mathbb{G}_{m,\widetilde{S}}) \ar[r]^{\kappa} \ar[d]^{H^1_{\acute{e} t}(i^*)} & H^1_{\acute{e} t}(\widetilde{S}_{\infty}\times\overline{\mathbb{Q}},i^*\mathbb{G}_{m,\widetilde{S}}) \ar[d]^{H^1_{\acute{e} t}(i^*)} \\
H^1_{\acute{e}t}(\widetilde{S}_{\infty}\times\overline{\mathbb{Q}},\boldsymbol{\mu}_{l^n})
\ar[r]  & H^1_{\acute{e}t}(\widetilde{S}_{\infty}\times\overline{\mathbb{Q}},\mathbb{G}_{m})
\ar[r]^{\kappa} & H^1_{\acute{e}t}(\widetilde{S}_{\infty}\times\overline{\mathbb{Q}},\mathbb{G}_{m}).}$$ 
We note that the restriction to the boundary factorises into $H^1_{\acute{e} t}(i^*) \circ R$, i.e.
$(H^1_{\acute{e} t}(i^*) \circ R) (\widetilde{L}) =
\widetilde{L}|_{\widetilde{S}_{\infty} \times \overline{\mathbb{Q}}}.$ By Section \ref{clb} we know  
$\widetilde{L}|_{\widetilde{S}_{\infty} \times \overline{\mathbb{Q}}} \in \mathbb{G}_{m,\overline{\mathbb{Q}}} = \text{Pic}^0(\widetilde{S}_{\infty} \times \overline{\mathbb{Q}}) \subset H^1_{\acute{e} t}(\widetilde{S}_{\infty} \times \overline{\mathbb{Q}},\mathbb{G}_{m})$
and even more $\widetilde{L}|_{\widetilde{S}_{\infty} \times \overline{\mathbb{Q}}} = \varepsilon^{-2}$.
 This restriction class $\widetilde{L}|_{\widetilde{S}_{\infty} \times
   \overline{\mathbb{Q}}}$ is the obstruction of the triviality of the
 extension. In the case it would vanish, then the lifting
 $\widetilde{c_1(L)}$ would come from an element in $H^1_{\acute{e}
   t,c}(\widetilde{S}\times\overline{\mathbb{Q}},\mathbb{G}_m)$, i.e. it would
 be Galois invariant. 

To prepare the proof of Theorem \ref{MTladic} below,
we must investigate some diagram chases in the diagram of Lemma \ref{BD}. This
exhibits how the Galois action on $[H^2_{\CHE,l}(S,\mathbb{Q}(1))(-)]^{\vee}$
and the above Kummer sequence in the cohomology are linked. We start with $c_1^{(l^n)}(\widetilde{L}) \in
H^2_{\acute{e}
  t}(\widetilde{S}\times\overline{\mathbb{Q}},\boldsymbol{\mu}_{l^n})$ in the
big diagram. This 
comes from the line bundle $\widetilde{L} \in H^1_{\acute{e}
  t}(\widetilde{S}\times\overline{\mathbb{Q}},\mathbb{G}_m)$. If we send
$\widetilde{L}$ to $H^1_{\acute{e}
  t}(\widetilde{S}_{\infty}\times\overline{\mathbb{Q}},i^*\mathbb{G}_{m,\widetilde{S}})$, we get a map $\varrho$ from $H^2_{\acute{e} t}(\widetilde{S}\times\overline{\mathbb{Q}},\boldsymbol{\mu}_{l^n})$ to $H^1_{\acute{e} t}(\widetilde{S}_{\infty}\times\overline{\mathbb{Q}},i^*\mathbb{G}_{m,\widetilde{S}})$. But this is not precisely correct, as there is some arbitrariness caused by the liftings. To understand this we consider the map
$H^2_{\acute{e} t}(\alpha_1 \oplus \alpha_2): H^2_{\acute{e} t}(\widetilde{S}\times\overline{\mathbb{Q}},\boldsymbol{\mu}_{l^n}) \rightarrow H^2_{\acute{e} t}(\widetilde{S}\times\overline{\mathbb{Q}},\mathbb{G}_m) \oplus  H^2_{\acute{e} t}(\widetilde{S}_{\infty}\times\overline{\mathbb{Q}},\boldsymbol{\mu}_{l^n}),$ which is
induced by the right bottom corner of the big diagram.  The kernel $\text{Ker}(H^2_{\acute{e} t}(\alpha_1 \oplus \alpha_2))$ is the subgroup of those classes that are Chern classes of line bundles and that have furthermore trivial Chern class on the boundary $\widetilde{S}_{\infty}\times\overline{\mathbb{Q}}$.
Note that by construction 
$c_1^{(l^n)}(\widetilde{L}) \in \text{Ker}(H^2_{\acute{e} t}(\alpha_1 \oplus \alpha_2)).$
And secondly, consider the map $H^1_{\acute{e} t}(\beta)$ defined by
$R \circ \kappa : H^1_{\acute{e} t}(\widetilde{S}\times\overline{\mathbb{Q}},\mathbb{G}_m)
\rightarrow
H^1_{\acute{e}t}(\widetilde{S}_{\infty}\times\overline{\mathbb{Q}},i^*\mathbb{G}_{m,\widetilde{S}}).$
The notation $H^i_{\acute{e} t}(-)$ indicates that these morphisms actually
come from morphisms of sheaves in the above diagram. We want to construct two maps 
$\text{Ker}(H^2_{\acute{e} t}(\alpha_1 \oplus \alpha_2))\rightarrow
\text{Coker}(H^1_{\acute{e} t}(\beta)).$ Let us write the whole diagram again.
{\scriptsize $$\xymatrix@C=0.8em{ & & & 0 \ar[d] \ar[r] & H^1_{\acute{e} t}(\widetilde{S}\times\overline{\mathbb{Q}},\mathbb{G}_m) \ar[d] \\
 & & & H^1_{\acute{e} t}(\widetilde{S}_{\infty}\times\overline{\mathbb{Q}},\boldsymbol{\mu}_{l^n}) \ar[d] \ar[r] & H^1_{\acute{e} t}(\widetilde{S}_{\infty}\times\overline{\mathbb{Q}},i^*\mathbb{G}_{m,\widetilde{S}}) \ar[d] \\  
& H^1_{\acute{e} t,c}(S\times\overline{\mathbb{Q}},\mathbb{G}_m) \ar[d] \ar[r] & H^1_{\acute{e} t,c}(S\times\overline{\mathbb{Q}},\mathbb{G}_m) \ar[d] \ar[r] & H^2_{\acute{e} t,c}(S\times\overline{\mathbb{Q}},\boldsymbol{\mu}_{l^n}) \ar[d] \ar[r] & H^2_{\acute{e} t,c}(S\times\overline{\mathbb{Q}},\mathbb{G}_m) \ar[d] \\
0 \ar[r] \ar[d] & H^1_{\acute{e} t}(\widetilde{S}\times\overline{\mathbb{Q}},\mathbb{G}_m) \ar[d]_{R} \ar[r] & H^1_{\acute{e} t}(\widetilde{S}\times\overline{\mathbb{Q}},\mathbb{G}_m) \ar[d]_<<<<<<<{R} \ar[r] &  H^2_{\acute{e} t}(\widetilde{S} \times\overline{\mathbb{Q}},\boldsymbol{\mu}_{l^n}) \ar[d]^{H^2_{\acute{e} t}(\alpha_2)} \ar[r]^{H^2_{\acute{e} t}(\alpha_1)} & H^2_{\acute{e} t}(\widetilde{S}\times\overline{\mathbb{Q}},\mathbb{G}_m) \\
 H^1_{\acute{e} t}(\widetilde{S}_{\infty}\times\overline{\mathbb{Q}},\boldsymbol{\mu}_{l^n}) \ar[d] \ar[r] & H^1_{\acute{e} t}(\widetilde{S}_{\infty}\times\overline{\mathbb{Q}},i^*\mathbb{G}_{m,\widetilde{S}}) \ar[r]_{\kappa} \ar[dr]_{\overline{\kappa}} \ar[d] & H^1_{\acute{e} t}(\widetilde{S}_{\infty}\times\overline{\mathbb{Q}},i^*\mathbb{G}_{m,\widetilde{S}}) \ar@{->>}[d] \ar[r] & H^2_{\acute{e} t}(\widetilde{S}_{\infty}\times\overline{\mathbb{Q}},\boldsymbol{\mu}_{l^n}) \ar[d] &{\text{Ker}(H^2_{\acute{e} t}(\alpha_1 \oplus \alpha_2))} \ar@{_{(}->}[ul] \ar@/^3pc/@{.>}[uuu]^>>>>>>>>{\lambda} \ar@/_2.5pc/@{.>}[lll]_>>>>>>>>{\lambda} \ar@/^2.5pc/@{.>}[dll]^{\overline{\varrho}} \ar@/^2pc/@{.>}[ll]^(.7){\varrho} \\ 
H^2_{\acute{e} t,c}(S\times\overline{\mathbb{Q}},\boldsymbol{\mu}_{l^n}) \ar[r] & H^2_{\acute{e} t,c}(S\times\overline{\mathbb{Q}},\mathbb{G}_m) &{\text{Coker}(H^1_{\acute{e} t}(\beta))} & 0 & }$$}
Since we divide out the image of $H^1_{\acute{e} t}(\beta)$ we know that
$\overline{\varrho}: \text{Ker}(H^2_{\acute{e} t}(\alpha_1 \oplus \alpha_2))
\rightarrow \text{Coker}(H^1_{\acute{e} t}(\beta))$ is a well-defined map. But in the diagram we have the dotted arrows $\lambda$ from
$\text{Ker}(H^2_{\acute{e} t}(\alpha_1 \oplus \alpha_2))$ to $H^1_{\acute{e}
  t}(\widetilde{S}_{\infty}\times\overline{\mathbb{Q}},i^*\mathbb{G}_{m,\widetilde{S}})$, where that sits now in the right top corner. If we compose this with the map $\overline{\kappa}$, induced by the Kummer map, we get a morphism $\overline{\kappa} \circ \lambda$, which lands in the same quotient $\text{Coker}(H^1_{\acute{e} t}(\beta))$ of $H^1_{\acute{e} t}(\widetilde{S}_{\infty}\times\overline{\mathbb{Q}},i^*\mathbb{G}_{m,\widetilde{S}})$ as the $\overline{\varrho}$. The following lemma shows that these are equal.
\begin{Lem}\label{HA}

With  the above notations, the two maps $\overline{\kappa} \circ \lambda$ and $\overline{\varrho}$ are equal.
\end{Lem}
\begin{proof}
First we have to see that $\overline{\kappa} \circ \lambda$ is indeed well-defined.

For this regard again the above diagram of cohomology groups. We know that
$c_1^{(l^n)}(\widetilde{L}) \in \text{Ker}(H^2_{\acute{e} t}(\alpha_1 \oplus
\alpha_2)) \subset H^2_{\acute{e}
  t}(\widetilde{S}\times\overline{\mathbb{Q}},\boldsymbol{\mu}_{l^n})$, or any
other element, comes from a lifting in $H^2_{\acute{e}
  t,c}(S\times\overline{\mathbb{Q}},\boldsymbol{\mu}_{l^n})$. If we send this
lifting via the natural coboundary morphism to $H^2_{\acute{e}
  t,c}(S\times\overline{\mathbb{Q}},\mathbb{G}_m)$, the image cannot
vanish. Otherwise the lifting would be the class, coming from an element in
$H^1_{\acute{e} t,c}(S\times\overline{\mathbb{Q}},\mathbb{G}_m)$, and this
would contradict by exactness the non-vanishing of the boundary class of
$\widetilde{L}$ in $H^1_{\acute{e}
  t}(\widetilde{S}_{\infty}\times\overline{\mathbb{Q}},i^*\mathbb{G}_{m,\widetilde{S}})$. But we know that we can lift this element in $H^2_{\acute{e} t,c}(S\times\overline{\mathbb{Q}},\mathbb{G}_m)$ to $H^1_{\acute{e} t}(\widetilde{S}_{\infty}\times\overline{\mathbb{Q}},i^*\mathbb{G}_{m,\widetilde{S}})$, since $c_1^{(l^n)}(\widetilde{L}) \in \text{Ker}(H^2_{\acute{e} t}(\alpha_1)).$ So we are in the right top corner. Apply the Kummer map $\overline{\kappa}$ and we are done. We have to take into considerations the non-uniqueness of the liftings. The first step does not effect anything, since the lift to $H^2_{\acute{e} t,c}(S\times\overline{\mathbb{Q}},\boldsymbol{\mu}_{l^n})$ is unique up to $H^1_{\acute{e} t}(\widetilde{S}_{\infty}\times\overline{\mathbb{Q}},\boldsymbol{\mu}_{l^n})$ and this is the kernel of the Kummer map. The second lift to $H^1_{\acute{e} t}(\widetilde{S}_{\infty}\times\overline{\mathbb{Q}},i^*\mathbb{G}_{m,\widetilde{S}})$ is unique up to elements in $H^1_{\acute{e} t}(\widetilde{S}\times\overline{\mathbb{Q}},\mathbb{G}_m)$. Here we use that we actually want to land into $\text{Coker}(H^1_{\acute{e} t}(\beta))$, hence this ambiguity is divided out. 
 Now I show that $\overline{\kappa} \circ \lambda$ and the
$\overline{\varrho}$ have to coincide. Since the rows in the diagram are
exact, it is sufficient to consider the map $\text{Ker}(H^2_{\acute{e}
  t}(\alpha')) \rightarrow \text{Coker}(H^1_{\acute{e} t}(\beta)),$ where 
$\text{Ker}(H^2_{\acute{e} t}(\alpha'))=\text{Ker}\left(H^2_{\acute{e}
    t,c}(S\times\overline{\mathbb{Q}},\boldsymbol{\mu}_{l^n}) \rightarrow
  H^2_{\acute{e}
    t}(\widetilde{S}\times\overline{\mathbb{Q}},\mathbb{G}_m)\right).$ 
Then we have a ``snake diagram''
{\scriptsize $$\xymatrix@C=2em{ & & & \text{Ker}(H^2_{\acute{e} t}(\alpha')) \ar[d] & \\
 & H^1_{\acute{e} t}(\widetilde{S}\times\overline{\mathbb{Q}},\mathbb{G}_m)
 \ar[d]^{H^1_{\acute{e} t}(\beta)} \ar[r] & H^1_{\acute{e}
 t}(\widetilde{S}\times\overline{\mathbb{Q}},\mathcal{C}(\alpha')) \ar[d]
 \ar[r] \ar[dl]^{H^1_{\acute{e} t}(u)} & H^2_{\acute{e} t,c}(S\times\overline{\mathbb{Q}},\boldsymbol{\mu}_{l^n}) \ar[d]^{H^2_{\acute{e} t}(\alpha')} \ar[r]^{H^2_{\acute{e} t}(\alpha')} \ar[dl]^{H^2_{\acute{e} t}(v)} & H^2_{\acute{e} t}(\widetilde{S}\times\overline{\mathbb{Q}},\mathbb{G}_m) \\
H^1_{\acute{e} t}(\widetilde{S}\times\overline{\mathbb{Q}},\mathbb{G}_m) \ar[r]^>>>>>{H^1_{\acute{e} t}(\beta)} & H^1_{\acute{e} t}(\widetilde{S}_{\infty}\times\overline{\mathbb{Q}},i^*\mathbb{G}_{m,\widetilde{S}}) \ar[r] \ar[d]  & H^2_{\acute{e} t}(\widetilde{S}\times\overline{\mathbb{Q}},\mathcal{K}(\beta)) \ar[r] & H^2_{\acute{e} t}(\widetilde{S}\times\overline{\mathbb{Q}},\mathbb{G}_m) & \\
 & \text{Coker}(H^1_{\acute{e} t}(\beta)) & & & } $$}
which is induced by the following exact diagram of sheaves on $\widetilde{S} \times\overline{\mathbb{Q}}$,
$$\xymatrix{0 \ar[r]  & j_!j^*\boldsymbol{\mu}_{l^n,\widetilde{S}}
  \ar[r]^<<<<<{\alpha'} \ar[d]^{u} &  {\mathbb{G}_{m,\widetilde{S}}} \ar[r] \ar@{=}[d] &
  {\mathcal{C}(\alpha')} \ar[r] \ar[d]^{v} & 0 \\
0 \ar[r] &  {\mathcal{K}(\beta)} \ar[r] &  {\mathbb{G}_{m,\widetilde{S}}} \ar[r]^<<<<{\beta} & {i_*i^*\mathbb{G}_{m,\widetilde{S}}} \ar[r] & 0.}$$
With the notation ${\mathcal{C}(\alpha')}$ for the cokernel of $\alpha'$, and
  ${\mathcal{K}(\beta)}$ for the kernel of $\beta$. Now one just imitate the
  proof of the snake lemma to get the well-defined map $\text{Ker}(H^2_{\acute{e}
  t}(\alpha')) \rightarrow \text{Coker}(H^1_{\acute{e} t}(\beta))$. The diagram does not
  fulfill exactly the assumption of the snake lemma, but here we are in a
  special case, which makes the things work. The ``snake diagram'' explains to
  us the two maps from $\text{Ker}(H^2_{\acute{e} t}(\alpha_1 \oplus
  \alpha_2))$ to $\text{Coker}(H^1_{\acute{e} t}(\beta))$ via the left bottom
  corner and the right top corner of the big diagram. First of all, this is due to the fact that there are two morphisms 
$H^1_{\acute{e}t}(\widetilde{S}\times\overline{\mathbb{Q}},\mathcal{C}(\alpha'))
  \rightarrow H^2_{\acute{e}
  t}(\widetilde{S}\times\overline{\mathbb{Q}},\mathcal{K}(\beta))$. One factorisation is given if we start with $H^1_{\acute{e} t}(\widetilde{S}\times\overline{\mathbb{Q}},\mathcal{C}(\alpha'))$ and go via the natural coboundary map to
$H^2_{\acute{e} t}(S\times\overline{\mathbb{Q}},j_!j^*\boldsymbol{\mu}_{l^n,\widetilde{S}}) = H^2_{\acute{e} t,c}(S\times\overline{\mathbb{Q}},\boldsymbol{\mu}_{l^n})$
 and compose this with the morphism induced on the second cohomology. And for
 the other round we start with $H^1_{\acute{e}
   t}(\widetilde{S}\times\overline{\mathbb{Q}},\mathcal{C}(\alpha'))$ and go
 via the induced map to $H^1_{\acute{e}
   t}(\widetilde{S}\times\overline{\mathbb{Q}},i_*i^*\mathbb{G}_{m,\widetilde{S}})$ and apply the other coboundary map to $H^2_{\acute{e} t}(\widetilde{S}\times\overline{\mathbb{Q}},\mathcal{K}(\beta))$. Now we are left to show that we get $\overline{\varrho}$ in the first case and $\overline{\kappa} \circ \lambda$ in the second one, so eventually by the commutativity they have to be equal.
To do so we observe that $u$ and $v$ in the above diagram can be factorise
in two ways: for $\varrho$ we have 
$$\xymatrix{0 \ar[r]  & j_!j^*\boldsymbol{\mu}_{l^n,\widetilde{S}}\ar@/_2pc/[dd]_<<<<{u}
  \ar[r]^<<<<<{\alpha'} \ar[d] &  {\mathbb{G}_{m,\widetilde{S}}} \ar[r] \ar@{=}[d] &
  {\mathcal{C}(\alpha')} \ar[r] \ar[d] \ar@/^2pc/[dd]^<<<<{v} & 0 \\
0 \ar[r] & {\boldsymbol{\mu}_{l^n,\widetilde{S}}} \ar[r]^{\alpha_1} \ar[d] & {\mathbb{G}_{m,\widetilde{S}}} \ar@{=}[d] \ar[r] & {\mathbb{G}_{m,\widetilde{S}}} \ar[d] \ar[r] & 0 \\
0 \ar[r] &  {\mathcal{K}(\beta)} \ar[r] &  {\mathbb{G}_{m,\widetilde{S}}}
  \ar[r]^<<<<{\beta} & {i_*i^*\mathbb{G}_{m,\widetilde{S}}} \ar[r] & 0.}$$
and for $\lambda$
$$\xymatrix{0 \ar[r]  & j_!j^*\boldsymbol{\mu}_{l^n,\widetilde{S}}\ar@/_2pc/[dd]_<<<<{u}
  \ar[r]^<<<<<{\alpha'} \ar[d] &  {\mathbb{G}_{m,\widetilde{S}}} \ar[r] \ar@{=}[d] &
  {\mathcal{C}(\alpha')} \ar[r] \ar[d] \ar@/^2pc/[dd]^<<<<{v} & 0 \\
0 \ar[r] &  {j_!j^* \mathbb{G}_{m,\widetilde{S}}} \ar[r] \ar[d] & {\mathbb{G}_{m,\widetilde{S}}} \ar@{=}[d] \ar[r] & {i_*i^*\mathbb{G}_{m,\widetilde{S}}} \ar[d] \ar[r] & 0 \\
0 \ar[r] &  {\mathcal{K}(\beta)} \ar[r] &  {\mathbb{G}_{m,\widetilde{S}}}
  \ar[r]^<<<<{\beta} & {i_*i^*\mathbb{G}_{m,\widetilde{S}}} \ar[r] & 0.}$$
This is just due to the commutativity of our initial diagram of sheaves.
\end{proof}
We come back to our situation with the
\begin{Cor}\label{i*e}

Let $H^1_{\acute{e} t}(i^*)$ be as above. Let $c_1^{(l^n)}(\widetilde{L}) \in
\text{Ker}(H^2_{\acute{e} t}(\alpha_1 \oplus \alpha_2))$ and let $\varepsilon
\in \mathcal{O}_F^*$ be as fixed in the very beginning. Then $ H^1_{\acute{e} t}(i^*)\left((\kappa \circ \lambda) (c_1^{(l^n)}(\widetilde{L}))\right) = H^1_{\acute{e} t}(i^*)\left(\varrho(c_1^{(l^n)}(\widetilde{L}))\right) = \varepsilon^{-2}. $ 
\end{Cor}
\begin{proof}
First recall the construction of $\lambda$ (Lemma \ref{HA}), i.e. we have the following piece in the big diagram
 $$\xymatrix{
H^1_{\acute{e} t}(\widetilde{S}\times\overline{\mathbb{Q}},\mathbb{G}_m) \ar[d] \ar[r] \ar[dr]^{H^1_{\acute{e} t}(\beta)} & H^1_{\acute{e} t}(\widetilde{S}\times\overline{\mathbb{Q}},\mathbb{G}_m) \ar[d] \\
H^1_{\acute{e} t}(\widetilde{S}_{\infty}\times\overline{\mathbb{Q}},i^*\mathbb{G}_{m,\widetilde{S}}) \ar[r]_{\kappa} \ar[d] & H^1_{\acute{e} t}(\widetilde{S}_{\infty}\times\overline{\mathbb{Q}},i^*\mathbb{G}_{m,\widetilde{S}}) \\
 H^2_{\acute{e} t,c}(S\times\overline{\mathbb{Q}},\mathbb{G}_m) & }$$
The lift from $H^2_{\acute{e} t,c}(S\times\overline{\mathbb{Q}},\mathbb{G}_m)$
to $H^1_{\acute{e}
  t}(\widetilde{S}_{\infty}\times\overline{\mathbb{Q}},i^*\mathbb{G}_{m,\widetilde{S}})$ is well-defined up to an image coming from $H^1_{\acute{e} t}(\widetilde{S}\times\overline{\mathbb{Q}},\mathbb{G}_m)$, but now we choose this such that the lift goes under $\kappa$ to $\varrho(c_1^{(l^n)}(\widetilde{L}))$. We can do this, since we know that the image of the lift via $\kappa$ is also only well-defined up to the image of $H^1_{\acute{e} t}(\beta)$. We know that $\varrho(c_1^{(l^n)}(\widetilde{L})) \in H^1_{\acute{e} t}(\widetilde{S}_{\infty} \times \overline{\mathbb{Q}},i^*\mathbb{G}_{m,\widetilde{S}})$ comes from the line bundle $\widetilde{L} \in H^1_{\acute{e} t}(\widetilde{S}\times\overline{\mathbb{Q}},\mathbb{G}_m)$.
Since the restriction map factorises into $H^1_{\acute{e} t}(i^*) \circ R$, we have that  
$H^1_{\acute{e} t}(i^*)\left(\varrho(c_1^{(l^n)}(\widetilde{L}))\right)$
equals the restriction of $\widetilde{L}$ to the boundary
$\widetilde{S}_{\infty} \times \overline{\mathbb{Q}}$, and we computed this in Lemma \ref{LEMCLB}, i.e.
$H^1_{\acute{e} t}(i^*)\left(\varrho(c_1^{(l^n)}(\widetilde{L}))\right) =
\widetilde{L}|_{\widetilde{S}_{\infty} \times \overline{\mathbb{Q}}} = \left( \widetilde{L_1}^{-1} \otimes \widetilde{L_2} \right)|_{\widetilde{S}_{\infty} \times \overline{\mathbb{Q}}} = \varepsilon^{-2}. $
\end{proof}
Let us now turn over to the Galois action and our main theorem. Recall that our big goal is to determine the Galois representations 
$ \sigma \mapsto \left(\begin{smallmatrix} \chi_D(\sigma) & * \\ 0 & \alpha^{-1}(\sigma) \end{smallmatrix}\right),$
which come from the $l$-adic realisations of our Kummer-Chern-Eisenstein
motive $[H^2_{\CHE}(S,\mathbb{Q}(1))(-)]$. The above discussion gives the star
$*$, and shows that it comes indeed from a Kummer extension of our real
quadratic field $F$. We sum this up in the 
\begin{Thm}\label{MTladic}

Let
$\varepsilon \in \mathcal{O}_F^*$ be as fixed in the very beginning and define
$\widetilde{\varepsilon}:= \varepsilon^{-\frac{1}{2}\zeta_F(-1)^{-1}}$. Then
$[H^2_{\CHE,l}(S,\mathbb{Q}(1))(-)]$ is $\widetilde{\varepsilon}$. The corresponding $l$-adic Galois representation comes from the Kummer field extension $F\left(\sqrt[l^{\infty}]{\widetilde{\varepsilon}}, \zeta_{l^{\infty}}\right)/F$ of $F$ attached to $\widetilde{\varepsilon}$, i.e.
$\sigma \mapsto
\left(\begin{smallmatrix} \chi_D(\sigma) & \tau_{\widetilde{\varepsilon}}(\sigma) \alpha^{-1}(\sigma) \\
0 & \alpha^{-1}(\sigma) \end{smallmatrix}\right),$ 
where $\tau_{\epsilon}(\sigma)$ is defined by 
$ \frac{\sigma\left(\sqrt[l^{\infty}]{\widetilde{\varepsilon}}\right)}{\sqrt[l^{\infty}]{\widetilde{\varepsilon}}} = \zeta_{l^{\infty}}^{\tau_{\widetilde{\varepsilon}}(\sigma)}.$
\end{Thm}
\begin{proof}

By Kummer theory we know that the first assertion follows from the second
one. 

We take off with the cocycle
$\sigma \left(\widetilde{c_1(L)}\right) -  \widetilde{c_1(L)},$
where the class $\widetilde{c_1(L)}$ is the lifting of $c_1(\widetilde{L}) \in H^2_{\acute{e} t}(\widetilde{S}\times\overline{\mathbb{Q}},\mathbb{Q}_l(1))$ (compare Lemma \ref{LEMLL}).
That cocycle gives the extension class of the dual Galois module
$[H^2_{\CHE,\acute{e} t}(S,\mathbb{Q}_l(1))(-)]^{\vee}$. We know that the dual motive sits in the following sequence (see Lemma \ref{MOTKK})
$$0  \rightarrow {
  \underset{n}{\varprojlim}H^1_{\acute{e}t}(\widetilde{S}_{\infty}\times\overline{\mathbb{Q}},\boldsymbol{\mu}_{l^n})
  \otimes \mathbb{Q}_l } \rightarrow  { \underset{n}{\varprojlim}
  H^2_{\acute{e} t,c}(S\times\overline{\mathbb{Q}},\boldsymbol{\mu}_{l^n})
  \otimes \mathbb{Q}_l} \rightarrow  {\text{Im}(f_1)} \rightarrow 0,$$ 
where this sequence itself is part of the big diagram of Lemma \ref{BD}.
 Now we apply the diagram chases in this diagram (loc. cit.). For start, we go
 again back to finite coefficients, i.e. consider the class
 $c_1^{(l^n)}(\widetilde{L}) \in H^2_{\acute{e}
   t}(\widetilde{S}\times\overline{\mathbb{Q}},\boldsymbol{\mu}_{l^n})$, and
 we write the action now multiplicatively. Furthermore, we should remind
 ourselves that $L$ is a certain power of $L$ - see Remark \ref{RMKPOT}. For
 the sake of simplicity, write in the following for this power again $L$. So we get the element ${\sigma}\left(\widetilde{c_1^{(l^n)}(L)}\right) \left(\widetilde{c_1^{(l^n)}(L)}\right)^{-1}$ in $H^2_{\acute{e} t,c}(S\times\overline{\mathbb{Q}},\boldsymbol{\mu}_{l^n})$, which vanishes under the restriction map to $H^2_{\acute{e} t}(\widetilde{S}\times\overline{\mathbb{Q}},\boldsymbol{\mu}_{l^n})$, i.e. it is in the kernel, and this is $H^1_{\acute{e} t}(\widetilde{S}_{\infty}\times\overline{\mathbb{Q}},\boldsymbol{\mu}_{l^n})$. By construction of $\lambda$, more precisely by the commutativity of 
$$\xymatrix{ H^1_{\acute{e} t}(\widetilde{S}_{\infty}\times\overline{\mathbb{Q}},\boldsymbol{\mu}_{l^n}) \ar[d] \ar[r] & H^1_{\acute{e} t}(\widetilde{S}_{\infty}\times\overline{\mathbb{Q}},i^*\mathbb{G}_{m,\widetilde{S}}) \ar[d] \\ H^2_{\acute{e} t,c}(S\times\overline{\mathbb{Q}},\boldsymbol{\mu}_{l^n}) \ar[r] & H^2_{\acute{e} t,c}(S\times\overline{\mathbb{Q}},\mathbb{G}_m)}$$
we conclude furthermore that up to an $l^n$-torsion element
$t^{\text{(obs)}}_{l^n}$ that we have the equality 
$ {\sigma}\left(\widetilde{c_1^{(l^n)}(L)}\right)
\left(\widetilde{c_1^{(l^n)}(L)}\right)^{-1} = \sigma \left( \lambda
  (c_1^{(l^n)}(\widetilde{L})) \right) \left( \lambda
  (c_1^{(l^n)}(\widetilde{L})) \right)^{-1}.$ 
The ambiguity given by $t^{\text{(obs)}}_{l^n}$ comes from the line bundles in
$H^1_{\acute{e} t}(\widetilde{S} \times\overline{\mathbb{Q}},\mathbb{G}_{m})$
that map via $R$ to $\sigma \left( \lambda (c_1^{(l^n)}(\widetilde{L}))
\right) \left( \lambda (c_1^{(l^n)}(\widetilde{L})) \right)^{-1}$. I call this
torsion element $t^{\text{(obs)}}_{l^n}$ with (``obs'' :=
``obstruction''). These line bundles, causing the trouble are also $l^n$ -
torsion elements in $H^1_{\acute{e} t}(\widetilde{S}
\times\overline{\mathbb{Q}},\mathbb{G}_{m})$. Now $H^1_{\acute{e}
  t}(\widetilde{S} \times\overline{\mathbb{Q}},\mathbb{G}_{m})$ is of finite
rank. (Note that $\Pic^0(\widetilde{S}\times\overline{\mathbb{Q}})$ vanishes.) Thus we
can only have finitely many $t^{\text{(obs)}}_{l^n}$'s in the image. Assume
that these line bundles are all of $l^k$-torsion. Otherwise, we choose a
suitable power, i.e. take $k \gg 0$. So let us consider again our situation with $l^{n}$-coefficients for $n \geq k$, 
$${\sigma}\left(\widetilde{c_1^{(l^{n})}(L)}\right) \left(\widetilde{c_1^{(l^{n})}(L)}\right)^{-1} = \sigma \left( \lambda (c_1^{(l^{n})}(\widetilde{L})) \right) \left( \lambda (c_1^{(l^{n})}(\widetilde{L})) \right)^{-1} \cdot t^{\text{(obs)}}_{l^{k}}.$$
Now if we raise this to $l^k$-th power, the obstruction vanishes, i.e.
$${\sigma}\left(\widetilde{c_1^{(l^{n-k})}(L)}\right) \left(\widetilde{c_1^{(l^{n-k})}(L)}\right)^{-1} = \sigma \left( \lambda (c_1^{(l^{n-k})}(\widetilde{L})) \right) \left( \lambda (c_1^{(l^{n-k})}(\widetilde{L})) \right)^{-1}.$$  
Hence in the limit 
${\sigma}\left(\widetilde{c_1^{(l^{\infty})}(L)}\right)
\left(\widetilde{c_1^{(l^{\infty})}(L)}\right)^{-1} = \sigma \left( \lambda
  (c_1^{(l^{\infty})}(\widetilde{L})) \right) \left( \lambda
  (c_1^{(l^{\infty})}(\widetilde{L})) \right)^{-1},$ 
i.e.
${\sigma}\left(\widetilde{c_1(L)}\right) \left(\widetilde{c_1(L)}\right)^{-1}
= \sigma \left( \lambda (c_1(\widetilde{L})) \right) \left( \lambda
  (c_1(\widetilde{L})) \right)^{-1}.$ By Corollary \ref{i*e} we know $(\kappa \circ \lambda) (c_1^{(l^n)}(\widetilde{L})) = \varrho(c_1^{(l^n)}(\widetilde{L}))$ and therefore
$\sigma \left( \lambda (c_1(\widetilde{L})) \right) \left( \lambda (c_1(\widetilde{L})) \right)^{-1} = \sigma\left(\sqrt[l^{\infty}]{\varrho(c_1(\widetilde{L}))}\right)\left(\sqrt[l^{\infty}]{\varrho(c_1(\widetilde{L}))}\right)^{-1},$
where $\sqrt[l^{\infty}]{-}$ indicates the preimage $\kappa^{-1}$ of the
Kummer map. We conclude that we have the equality 
${\sigma}\left(\widetilde{c_1(L)}\right) \left(\widetilde{c_1(L)}\right)^{-1}
=
\sigma\left(\sqrt[l^{\infty}]{\varrho(c_1(\widetilde{L}))}\right)\left(\sqrt[l^{\infty}]{\varrho(c_1(\widetilde{L}))}\right)^{-1}$ 
in $\underset{n}{\varprojlim}H^1_{\acute{e} t}(\widetilde{S}_{\infty}\times\overline{\mathbb{Q}},\boldsymbol{\mu}_{l^n}) \otimes \mathbb{Q}_l$. This group is the bottom of our extension $H^1_{\acute{e} t}(\widetilde{S}\times\overline{\mathbb{Q}},\mathbf{R}j_*\mathbb{Q}_l/j_!\mathbb{Q}_l) \otimes \mathbb{Q}_l(1) \simeq \mathbb{Q}_l(1)$.

Now we have to determine this element in $\mathbb{Q}_l(1)$. For this consider again the diagram (Section \ref{liftings})
$$\xymatrix{ & H^1_{\acute{e} t}(\widetilde{S}\times\overline{\mathbb{Q}},\mathbb{G}_m) \ar[d]^{R} \ar[r] & H^1_{\acute{e} t}(\widetilde{S}\times\overline{\mathbb{Q}},\mathbb{G}_m) \ar[d]^{R} \\
H^1_{\acute{e} t}(\widetilde{S}_{\infty}\times\overline{\mathbb{Q}},\boldsymbol{\mu}_{l^n}) \ar[r] \ar@{=}[d] & H^1_{\acute{e} t}(\widetilde{S}_{\infty}\times\overline{\mathbb{Q}},i^*\mathbb{G}_{m,\widetilde{S}}) \ar[r]^{\kappa} \ar[d]^{H^1_{\acute{e} t}(i^*)} & H^1_{\acute{e} t}(\widetilde{S}_{\infty}\times\overline{\mathbb{Q}},i^*\mathbb{G}_{m,\widetilde{S}}) \ar[d]^{H^1_{\acute{e} t}(i^*)} \\
H^1_{\acute{e} t}(\widetilde{S}_{\infty}\times\overline{\mathbb{Q}},\boldsymbol{\mu}_{l^n}) \ar[r] & H^1_{\acute{e} t}(\widetilde{S}_{\infty}\times\overline{\mathbb{Q}},\mathbb{G}_{m}) \ar[r]^{\kappa} & H^1_{\acute{e} t}(\widetilde{S}_{\infty}\times\overline{\mathbb{Q}},\mathbb{G}_{m}) \\
 & {\mathbb{G}_m} \ar@{^{(}->}[u] \ar[r] & {\mathbb{G}_m} \ar@{^{(}->}[u].}$$ 
If we apply the above Corollary \ref{i*e}, we get 
\begin{eqnarray*} 
H^1_{\acute{e} t}(i^*)\left(\sigma\left(\sqrt[l^n]{\varrho(c_1^{(l^n)}(\widetilde{L}))}\right)\left(\sqrt[l^n]{\varrho(c_1^{(l^n)}(\widetilde{L}))}\right)^{-1}\right) & = & \sigma\left(\sqrt[l^n]{\varepsilon^{-2}}\right)\left(\sqrt[l^n]{\varepsilon^{-2}}\right)^{-1},
\end{eqnarray*}
where 
$\sigma\left(\sqrt[l^n]{\varepsilon^{-2}}\right)\left(\sqrt[l^n]{\varepsilon^{-2}}\right)^{-1} \in F\left(\sqrt[l^n]{\varepsilon^{-2}}, \zeta_{l^n} \right)$ 
is an element in the Kummer extension of degree $l^n$ attached to
$\varepsilon^{-2}$. Moreover, this is exactly the Galois cocycle attached to
this field extension. Now the two Kummer maps above have the same kernel,
which is 
$H^1_{\acute{e}
  t}(\widetilde{S}_{\infty}\times\overline{\mathbb{Q}},\boldsymbol{\mu}_{l^n})$, and hence we see that in the limit $\underset{n}{\varprojlim}$ that
$
{\sigma}\left(\widetilde{c_1(L)}\right) \left(\widetilde{c_1(L)}\right)^{-1} = \sigma\left(\sqrt[l^{\infty}]{\varrho(c_1(\widetilde{L}))}\right)\left(\sqrt[l^{\infty}]{\varrho(c_1(\widetilde{L}))}\right)^{-1} = \sigma\left(\sqrt[l^{\infty}]{\varepsilon^{-2}}\right)\left(\sqrt[l^{\infty}]{\varepsilon^{-2}}\right)^{-1}.
$
All this yields eventually to the Galois representation 
$ \sigma  \mapsto
\left(\begin{smallmatrix} \alpha(\sigma) &  \tau_{\varepsilon^{-2}}(\sigma)\alpha^{-1}(\sigma) \\
0 &  \chi_D(\sigma)  \end{smallmatrix}\right),$
where 
$\frac{\sigma\left(\sqrt[l^{\infty}]{\varepsilon^{-2}}\right)}{\sqrt[l^{\infty}]{\varepsilon^{-2}}} = \zeta_{l^{\infty}}^{ \tau_{\varepsilon^{-2}}(\sigma)}.$

We know that the Galois cocycle for the actual module
$[H^2_{\CHE,l}(S,\mathbb{Q}(1))(-)]$ in
$\Ext^1_{\mathcal{MGAL}}(\mathbb{Q}_l(-1),\mathbb{Q}_l(0)\chi_D)$ is  given by
$ {\sigma}\left(-\widetilde{c_1(L)}\right) -
\left(-\widetilde{c_1(L)}\right),$ 
i.e. we can play the same game with the inverse. But here we have to be
careful as we have to normalise the generator after dualising, see Remark \ref{RMKDM}, i.e. we have to multiply with
$-\frac{1}{4\zeta_F(-1)}$. Therefore we get
$(\varepsilon^2)^{-\frac{1}{4\zeta_F(-1)}} = \widetilde{\varepsilon}$ and our
Galois representation is
$ \sigma  \mapsto
\left(\begin{smallmatrix}  \chi_D(\sigma) & \tau_{\widetilde{\varepsilon}}(\sigma)\alpha^{-1}(\sigma)\\
0 & \alpha^{-1}(\sigma) \end{smallmatrix}\right).$
\end{proof}
\begin{Rmk}\label{RMKFQ}
Since the restriction of $L_1 \otimes L_2$ is trivial on
$\widetilde{S}_{\infty}\times\overline{\mathbb{Q}}$ (Lemma \ref{LEMCLB}) and
its Chern class generates the $(+1)$-eigenspace $\mathbb{Q}(0)$ (Remark
\ref{RMKGEN}), we get a three dimensional representation $\sigma \mapsto  \left(\begin{smallmatrix} 1 & 0 & 0\\ 0 & \chi_D(\sigma) &  \tau_{\widetilde{\varepsilon}}(\sigma) \alpha^{-1}(\sigma)\\ 0 & 0 & \alpha^{-1}(\sigma) \end{smallmatrix}\right),$
that is induced by 
$0 \rightarrow \mathbb{Q}_l(0)\oplus \mathbb{Q}_l(0)\chi_D \rightarrow H^2_{\CHE,l}(S,\mathbb{Q}(1)) \rightarrow \mathbb{Q}_l(-1) \rightarrow 0.$
\end{Rmk}
Let us note that this is the realisation of the Kummer motive $K \langle \widetilde{\varepsilon} \rangle$ attached to our
$\widetilde{\varepsilon}$, i.e. $K \langle \widetilde{\varepsilon} \rangle_l = [H^2_{\CHE,l}(S,\mathbb{Q}(1))(-)] = \widetilde{\varepsilon}.$

\section{The Hodge-de Rham Realisation}\label{HdR}
In this chapter we compute the extension class 
$[H^2_{\CHE,\infty}(S,\mathbb{Q}(1))(-)]$ of the Hodge-de Rham realisation,
which is an element in $\Ext^1_{\mathcal{MH}d\mathcal{R}_{\mathbb{Q}}}(\mathbb{Q}(-1),\mathbb{Q}(0)\chi_D).$ 

In Section \ref{MMwC} we give the recipe that describes such an
element. We must understand the sections $s_B$ and $s_{dR}$ in our setting. 
Recall that Corollary \ref{CC} is the starting point of the construction of the motive, i.e. we have the sequence
$0\rightarrow H^2_!(S,\mathbb{Q}(1))\rightarrow H^2(S,\mathbb{Q}(1))\rightarrow \mathbb{Q}(-1)\rightarrow 0.$
And the top $\mathbb{Q}(-1)$ is given by the toroidal compactification
$\widetilde{S}$ as $H^2(\widetilde{S},\mathbf{R}j_*\mathbb{Q}/j_!\mathbb{Q}) \otimes
\mathbb{Q}(1).$ Now we consider the complex points $S(\mathbb{C})$ and let
$\partial\widetilde{S}_{\infty}(\mathbb{C})$ be the boundary of a suitable neighbourhood of $\widetilde{S}_{\infty}(\mathbb{C})$. 
By \cite{HLR}, 5, we have the exact sequence
$$ 0 \rightarrow H^1(\partial\widetilde{S}_{\infty}(\mathbb{C}),\mathbb{Q}) \rightarrow H^2_c(S(\mathbb{C}),\mathbb{Q}) \rightarrow H^2(S(\mathbb{C}),\mathbb{Q}) \rightarrow H^2(\partial\widetilde{S}_{\infty}(\mathbb{C}),\mathbb{Q})\rightarrow 0.$$
For the vanishing of $H^1(S(\mathbb{C}),\mathbb{Q})$, resp. $H^3_c(S(\mathbb{C}),\mathbb{Q})$  see Lemma \ref{LEMBLC}.
So the cokernel $H^2(\partial\widetilde{S}_{\infty}(\mathbb{C}),\mathbb{Q})$
is isomorphic to
$H^2(\widetilde{S}(\mathbb{C}),\mathbf{R}j_*\mathbb{Q}/j_!\mathbb{Q})$. By
\cite{AMRT}, I.5, or \cite{HLR}, 3 and 5, we know that
$\partial\widetilde{S}_{\infty}(\mathbb{C})$ is isomorphic to the boundary
$\partial S(\mathbb{C})^{BS}$ of the Borel-Serre compactification
$S(\mathbb{C}) \hookrightarrow S(\mathbb{C})^{BS}.$ This put us in the
position to describe the sections $s_B$ and $s_{dR}$ by Eisenstein cohomology,
and that is done in the next section.
\subsection{The Extension Class as Eisenstein Class} \label{Eisensteinclass}%
All this bases on Harder's notes \cite{?!}, see also \cite{SL2(O)} and
\cite{KAG}. Consider the short exact sequence
$$0 \rightarrow H^2_!(S(\mathbb{C}),\mathbb{C}) \rightarrow H^2(S(\mathbb{C})^{BS},\mathbb{C}) \rightarrow H^2(\partial S(\mathbb{C})^{BS},\mathbb{C}) \rightarrow 0,$$
which now comes from the Borel-Serre compactification $S(\mathbb{C}) \hookrightarrow S(\mathbb{C})^{BS}.$
The homotopy equivalence between $S(\mathbb{C})$ and $S(\mathbb{C})^{BS}$ induces the isomorphism between
$H^2(S(\mathbb{C}),\mathbb{C})$ and $H^2(S(\mathbb{C})^{BS},\mathbb{C}).$ For simplicity, we restrict the general setting concerning Eisenstein
cohomology to our situation, i.e. to the case of our group $G =
\Res_{F/\mathbb{Q}}(\GL_2/F)$ and the cohomology in degree two
$H^2(S(\mathbb{C}),\mathbb{C})$ with the constant coefficient system
$\mathbb{C}$. Let us start with the de Rham theorem (\cite{KAG}, Satz 3.7.8.2)
$ H^*(\partial S(\mathbb{C})^{BS},\mathbb{C}) \simeq
H^*(\mathfrak{g},K;\mathcal{C}^{\infty}(B(\mathbb{Q})\backslash
G(\mathbb{A}))),$ 
with the Lie algebra $\mathfrak{g} = \text{Lie}(G_{\infty})$ of $G_{\infty}$, the group $K=K_{\infty}$ and the standard Borel subgroup $B \subset G$, where the notation are as in
 Section \ref{HMS}. We define (as in \cite{KAG}, 5.2, or \cite{SL2(O)}, 2) 
$\Eis: \mathcal{C}^{\infty}(B(\mathbb{Q})\backslash G(\mathbb{A})) \rightarrow \mathcal{C}^{\infty}(G(\mathbb{Q})\backslash G(\mathbb{A}))$
by $\Eis(\psi):=\{g\mapsto\sum_{\underline{a}\in B(\mathbb{Q})\backslash
  G(\mathbb{Q})} \psi(\underline{a}g)\}.$ This induces a map
$H^2(\mathfrak{g},K;\mathcal{C}^{\infty}(B(\mathbb{Q})\backslash G(\mathbb{A}))) \rightarrow H^2(\mathfrak{g},K;\mathcal{C}^{\infty}(G(\mathbb{Q})\backslash G(\mathbb{A})))$ 
and even more a section
$H^2(\partial S(\mathbb{C})^{BS},\mathbb{C}) \rightarrow H^2(S(\mathbb{C}),\mathbb{C})$
of the restriction map to the boundary, see loc. cit. To describe this in more
detail, we go back to the Lie algebra cohomology. We have the isomorphism
(again by the de Rham theorem) 
$\Hom_K(\Lambda^2(\mathfrak{g}/\mathfrak{k}),\mathcal{C}^{\infty}(B(\mathbb{Q})\backslash
G(\mathbb{A}))) = \Omega^2(B(\mathbb{Q})\backslash
G(\mathbb{A})/K_{\infty}K_f) = \Omega^2(\partial S(\mathbb{C})^{BS}),$
where $\mathfrak{k} := \text{Lie}(K_{\infty})$ is the Lie algebra of
$K_{\infty}$. The product $G(\mathbb{R}) \simeq \GL_2(\mathbb{R}) \times \GL_2(\mathbb{R})$ induces a Hodge decomposition
$\Lambda^2(\mathfrak{g}/\mathfrak{k})  =
\Lambda^2(\mathfrak{g}_1/\mathfrak{k}_1) \oplus
(\Lambda^1(\mathfrak{g}_1/\mathfrak{k}_1) \otimes
\Lambda^1(\mathfrak{g}_2/\mathfrak{k}_2)) \oplus
\Lambda^2(\mathfrak{g}_2/\mathfrak{k}_2)$ for the exterior algebra $\Lambda^2(\mathfrak{g}/\mathfrak{k})$. Now we can choose a dual basis $\omega_{+,j}, \omega_{-,j}$ of $\Lambda^1(\mathfrak{g}_j/\mathfrak{k}_j)$, which corresponds to $dz_j,
d\overline{z}_j$ or $dx_j \pm idy_j$, $j=1,2$, see for example \cite{HLR}, 3,
or  \cite{HaLNM}, 4.3.3. 
 To simplify the notation, we denote 
$\delta z_j := \omega_{+,j}$ and $\delta \overline {z}_j := \omega_{-,j},$
this means for the cohomology classes 
$[dz_j] = [\delta z_j].$
Since 
$\{dz_1 \wedge d{z}_2, dz_1 \wedge d\overline{z}_2, d\overline{z}_1 \wedge
dz_2, d\overline{z}_1 \wedge d\overline{z}_2\}$ generate the same $\mathbb{C}$-vector space as 
$\{dx_1 \wedge dx_2, dx_1 \wedge dy_2, dx_2 \wedge dy_1, dy_1 \wedge dy_2\},$ 
we get in same manner elements
$\{\delta x_1\wedge \delta x_2, \delta x_1\wedge\delta y_2, \delta x_2\wedge \delta y_1, \delta y_1\wedge \delta y_2\}$
in
$\Hom_K(\Lambda^2(\mathfrak{g}/\mathfrak{k}),\mathcal{C}^{\infty}(B(\mathbb{Q})\backslash
G(\mathbb{A})))$, which form a dual basis. Now we have to consider the finite
part of the cohomology of the boundary $H^2(\partial
S(\mathbb{C})^{BS},\mathbb{C})$. Since we know that this is one dimensional,
the finite adelic part is just generated by a normed standard spherical
function $\psi_f$, i.e. $\psi_f(1) = 1$. More precisely, we have the
\begin{Lem}\label{LEMGEN}

The cohomology class $[(\delta x_1\wedge \delta x_2)\otimes \psi_f] \in H^2(\partial S(\mathbb{C})^{BS},\mathbb{C})$ of 
$(\delta x_1\wedge \delta x_2)\otimes \psi_f \in \Hom_K(\Lambda^2(\mathfrak{g}/\mathfrak{k}),\mathcal{C}^{\infty}(B(\mathbb{Q})\backslash G(\mathbb{A})))$ 
generates the cohomology $H^2(\partial S(\mathbb{C})^{BS},\mathbb{C})$ of the Borel-Serre boundary $\partial S(\mathbb{C})^{BS}$.  
\end{Lem}
\begin{proof}

See e.g. \cite{SL2(O)}, Proposition 1.1.
\end{proof}
Let us denote this generator $\omega_0$ by
$\omega_0 := \delta x_1\wedge \delta x_2.$ Then we have
\begin{Lem}\label{COR42}
Let $\omega=\delta x_1\wedge \delta x_2+\alpha_1 \cdot \delta x_1\wedge\delta y_2+\alpha_2 \cdot \delta x_2\wedge \delta y_1+ \beta \cdot \delta y_1\wedge \delta y_2$ 
be a form, which defines a closed form $\omega \otimes \psi_f \in
\Omega^2(\partial S(\mathbb{C})^{BS})$. Then its cohomology class $[\omega
\otimes \psi_f]$ in $H^2(\partial S(\mathbb{C})^{BS},\mathbb{C})$ equals
$[\omega_0 \otimes \psi_f]$, i.e. it is independent of the coefficients
$\alpha_1, \alpha_2$ and $\beta$. If we apply the Eisenstein operator to $\omega \otimes \psi_f$, we get a closed form
$\Eis(\omega \otimes \psi_f)\in \Hom_K(\Lambda^2(\mathfrak{g}/\mathfrak{k}),\mathcal{C}^{\infty}(G(\mathbb{Q})\backslash G(\mathbb{A}))),$
 and the class 
$[\Eis(\omega \otimes \psi_f)] \in H^2(S(\mathbb{C}),\mathbb{C})$ 
restricted to the boundary $\partial S(\mathbb{C})^{BS}$ is again 
$[\Eis(\omega \otimes \psi_f)]|_{\partial S(\mathbb{C})^{BS}} = [\omega \otimes \psi_f] \in H^2(\partial S(\mathbb{C})^{BS},\mathbb{C}),$ 
i.e. independent of the coefficients $\alpha_1, \alpha_2$ and $\beta$. And $\Eis$ is indeed a section for the restriction map.
\end{Lem}
\begin{proof}

See e.g. \cite{SL2(O)}, Theorem 2.1.
\end{proof} 
Let us compute the Eisenstein class.
\begin{Thm}\label{MT2}
Let
$\omega = \delta x_1\wedge \delta x_2+\alpha_1 \cdot \delta x_1\wedge\delta y_2+\alpha_2 \cdot \delta x_2\wedge \delta y_1+ \beta \cdot \delta y_1\wedge \delta y_2$  
be a closed two form $\omega\in\Omega^2(\partial S(\mathbb{C})^{BS})$. Then $[\Eis(\omega \otimes \psi_f)]\in H^2(S(\mathbb{C}),\mathbb{C})$ is
$[\Eis(\omega_0 \otimes \psi_f)] + \frac{h \sqrt{D} \cdot \log\varepsilon}{4
  \pi  \cdot \zeta_F(-1)} \left(\alpha_1 \cdot c_1(L_1) + \alpha_2 \cdot
  c_1(L_2)\right),$ 
where again $\varepsilon=\varepsilon_0^2 \in \mathcal{O}^*_F$ is our fixed
totally positive unit and $h$ the class number of $F$ (where we actually assume that $h=1$).
\end{Thm}
\begin{proof}

According to the last lemma above, the difference of the two sections $[\Eis(\omega \otimes \psi_f)]$ and $[\Eis(\omega_0 \otimes \psi_f)]$ vanishes under the restriction to the boundary, i.e.
$[\Eis(\omega \otimes \psi_f)]-[\Eis(\omega_0 \otimes \psi_f)] \in
H^2_!(S(\mathbb{C}),\mathbb{C}).$ 
In particular, it lies in the space $H^2_{\CH}(S(\mathbb{C}),\mathbb{C})$,
which is generated by the two Chern classes $c_1(L_1), c_1(L_2)$, see Section
\ref{Construction}. Hence we have to compute the relation
$[\Delta]:=[\Eis(\omega \otimes \psi_f)]-[\Eis(\omega_0 \otimes \psi_f)]=
\lambda_1 \cdot c_1(L_1) + \lambda_2 \cdot c_1(L_2).$ 
This is done in the same manner as in the proof of \cite{GL2}, Proposition
3.2.4 (see also \cite{SL2(O)}, 2). Recall that the $c_1(L_j)$ are the cohomology classes of 
$\zeta_j:=\frac{\delta x_j \wedge \delta y_j}{y_j^2},$
i.e.
$2 \pi  \cdot c_1(L_j) = [\zeta_j].$
By the use of loc. cit. (compare additionaly \cite{SL2(O)}, Proposition 2.3) we know that the $\zeta$'s are cohomologous to forms with compact support $\widetilde{\zeta_j}$, i.e. 
$\widetilde{\zeta_j} = \zeta_j - d\Psi_j,$
where 
$\Psi_j:=-f \cdot \frac{\delta x_j}{y_j},$
and where $f$ is a $\mathcal{C}^{\infty}$-function on the boundary that has support in the neighbourhood of the cusp, and is equal to one in a smaller neighbourhood (see loc. cit.). So $\Psi_j$ bounds $\zeta_j$ around the cusp. To get the coefficients $\lambda_1, \lambda_2$, we cup the above equation 
$[\Delta] = \lambda_1 \cdot c_1(L_1) + \lambda_2 \cdot c_1(L_2) = \frac{1}{2
  \pi } \left(\lambda_1 \cdot [\zeta_1] + \lambda_2 \cdot [\zeta_2]\right) = \lambda_1' \cdot [\zeta_1] + \lambda_2' \cdot [\zeta_2]$
with $[\widetilde{\zeta_2}]$. The cup product gives on the right hand side $ \left(\lambda'_1 [\zeta_1] +
   \lambda'_2 [\zeta_2]\right) \cup [\widetilde{\zeta_2}] = \lambda'_1
 [\zeta_1] \cup [\widetilde{\zeta_2}] + \left( \lambda'_2 [\zeta_2] \cup
   [\zeta_2] - \lambda'_2 [\zeta_2] \cup  d\Psi_2\right).$ Recall that here in our case $S_{K_0}(\mathbb{C}) = \Gamma \backslash (\mathfrak{H} \times \mathfrak{H}) ,$
with $\Gamma = \PSL_2(\mathcal{O}_F).$  
Then by loc. cit. the cup product reduces to  
$\lambda'_1 [\zeta_1] \cup [\zeta_2] =  \lambda'_1 \cdot \int_{\Gamma
  \backslash (\mathfrak{H} \times \mathfrak{H}) } \frac{\delta x_1 \wedge
  \delta y_1}{y_1^2} \wedge \frac{\delta x_2 \wedge \delta y_2}{y_2^2}  =
\lambda'_1 \cdot 4 \pi^2 \cdot \text{Vol}(\Gamma \backslash (\mathfrak{H}
\times \mathfrak{H})) =  \lambda'_1 \cdot 8 \pi^2 \cdot \zeta_F(-1),$
where the last equality is again Siegel's theorem. Now we have to compute the other side, which is a little bit more
delicate. For this we chop off the cusp at a certain level $c \gg 0$, i.e. we
consider the Borel-Serre compactification. Then the left hand side becomes 
$[\Delta] \cup [\widetilde{\zeta_2}] = \int_{\Gamma \backslash
  (\mathfrak{H}\times \mathfrak{H})_{\leq c}} \Delta \wedge
\widetilde{\zeta_2} = \int_{\Gamma \backslash (\mathfrak{H}\times
  \mathfrak{H})_{\leq c}} \Delta \wedge \left({\zeta_2} - d\Psi_2\right).$ And therefore 
$[\Delta] \cup [\widetilde{\zeta_2}] = - \int_{\Gamma \backslash
  (\mathfrak{H}\times \mathfrak{H})_{\leq c}} \Delta \wedge  d\Psi_2 = -
\int_{\partial \left(\Gamma \backslash (\mathfrak{H}\times \mathfrak{H})_{\leq
      c}\right)} \Delta \wedge  \Psi_2.$ According to \cite{SL2(O)}, 2, we know
$\Delta = (\alpha_1 \cdot \delta x_1\wedge\delta y_2+\alpha_2 \cdot \delta
x_2\wedge \delta y_1+ \beta \cdot \delta y_1\wedge \delta y_2) \otimes \psi_f+
O(\frac{1}{c^N}).$ We get
\begin{eqnarray*} 
- \int_{\partial \left(\Gamma \backslash (\mathfrak{H}\times \mathfrak{H})_{\leq c}\right)} \Delta \wedge \Psi_2  
& = & - \int_{\partial \left(\Gamma_{\infty} \backslash (\mathfrak{H}\times \mathfrak{H})_{\leq c}\right)} \Delta \wedge  \Psi_2 \\
 & = &  - \alpha_1 \int_{\partial \left(\Gamma_{\infty} \backslash (\mathfrak{H}\times \mathfrak{H})_{\leq c}\right)} \delta x_1 \wedge \delta x_2 \wedge \frac{\delta y_2}{y_2}.
\end{eqnarray*}
The second equality comes from the fact that we integrate over the boundary
circle, where the product $y_1 \cdot y_2$ of the imaginary parts is
constant. To calculate the latter integral, we recall (\cite{AMRT}, I.5) that
the boundary is a torus bundle over $\mathbb{S}^1$ with fibres isomorphic to
$\mathcal{O}_F \backslash \mathbb{R}^2$. Moreover, the base $\mathbb{S}^1$ is
given by the action of the units $\mathcal{O}_F^*$. According to our
orientation, we have to integrate in the (second) coordinate $y_2$ from $1$ to
$\varepsilon^{-1}$ with fibres $\mathcal{O}_F \backslash \mathbb{R}^2$,
i.e. the latter integral becomes 
$ \int_{\mathcal{O}_F \backslash \mathbb{R}^2} \delta x_1 \wedge \delta x_2
\cdot \int_1^{\varepsilon^{-1}} \frac{\delta y_2}{y_2} = - \sqrt{D}\cdot
\log\varepsilon,$ 
where the factor $\sqrt{D}$ is the Euclidean volume of our real quadratic
field $F$-note that $D \equiv 1 (\text{mod} \: 4)$. Therefore we get eventually  
$ [\Delta] \cup [\widetilde{\zeta_2}] =  - \int_{\partial
  \left(\Gamma_{\infty} \backslash (\mathfrak{H}\times \mathfrak{H})_{\leq
      c}\right)} \Delta \wedge \Psi_2 =  \alpha_1 \cdot \sqrt{D} \cdot
\log\varepsilon.$ Plugging in this into the above formula, we get
$ \lambda'_1 \cdot 8 \pi^2 \cdot \zeta_F(-1) = \alpha_1 \cdot \sqrt{D}\cdot \log \varepsilon.$ 
And we end up with 
$ \lambda'_1 =\frac{ \sqrt{D} \cdot \log\varepsilon}{ 8 \pi^2 \cdot
  \zeta_F(-1)} \cdot \alpha_1 .$ For the other coefficient $\lambda'_2$ we do
the same, but now we observe that we have to integrate from $1$ to
$\varepsilon$. Then the minus sign disappears, too. Altogether we end up with $[\Delta] =  \frac{  \sqrt{D} \cdot \log\varepsilon}{ 4 \pi
 \cdot  \zeta_F(-1)} \cdot \left(\alpha_1 \cdot c_1(L_1) + \alpha_2 \cdot
 c_1(L_2) \right)$. 
If we relax the assumption that $h=1$, we have to add up all the contributions from the different cusps.
\end{proof}
Now let us come back to the determination of the extension class 
$$[H^2_{\CHE,\infty}(S,\mathbb{Q}(1))(-)] \in \Ext^1_{\mathcal{MH}d\mathcal{R}_{\mathbb{Q}}}(\mathbb{Q}(-1),\mathbb{Q}(0)\chi_D).$$
In the above Theorem \ref{MT2} we have the term 
$  \frac{h \sqrt{D} \cdot \log\varepsilon}{4 \pi  \cdot \zeta_F(-1)} \left(\alpha_1 \cdot c_1(L_1) + \alpha_2  \cdot c_1(L_2)\right).$
This describes what happens, if we modify 
$\omega_0 \otimes \psi_f$ by a coboundary $d\phi$. Then $\omega_0 \otimes \psi_f$ and $(\omega_0 + d\phi) \otimes \psi_f$ have the same cohomology class 
$ [ (\omega_0 + d\phi) \otimes \psi_f ] = [\omega_0 \otimes \psi_f] \in
H^2(\partial S(\mathbb{C})^{BS},\mathbb{C}),$ 
but the image under Eisenstein may differ, and this is exactly given by the above term. 
This is the Hodge-de Rham extension  $[H^2_{\CHE,\infty}(S,\mathbb{Q}(1))(-)]$ in  $\Ext^1_{\mathcal{MH}d\mathcal{R}_{\mathbb{Q}}}(\mathbb{Q}(-1),\mathbb{Q}(0)\chi_D)$.
\begin{Thm}\label{HdREXT}
Let $\varepsilon \in \mathcal{O}_F^*$ be as fixed in the very beginning. Then the Hodge-de Rham realisation 
$[H^2_{\CHE,\infty}(S,\mathbb{Q}(1))(-)] \in
\Ext^1_{\mathcal{MH}d\mathcal{R}_{\mathbb{Q}}}(\mathbb{Q}(-1),\mathbb{Q}(0)\chi_D)$  
of our Kummer-Chern-Eisenstein motive is $- \frac{\log \varepsilon}{2 \cdot
  \zeta_F(-1)} = \log \widetilde{\varepsilon}.$ 
\end{Thm}
\begin{proof}
We have our diagram
$$\xymatrix{ 0 \ar[r] & {\mathbb{Q}(0)\chi_{D,B}} \otimes \mathbb{C} \ar[d]^{\simeq}_{\cdot (\sqrt{D})^{-1}} \ar[r] & H^2_{\CHE,B}(S,\mathbb{Q}(1))(-) \otimes \mathbb{C} \ar[d]^{\simeq}_{I_{\infty}} \ar[r] & {\mathbb{Q}(-1)_B} \otimes \mathbb{C} \ar[d]^{\simeq}_{\cdot (2\pi i)^{-1}} \ar[r] \ar@/_1pc/[l]_{s_B} & 0 \\
0 \ar[r] & {\mathbb{Q}(0)\chi_{D,dR}} \otimes \mathbb{C} \ar[r] & H^2_{\CHE,dR}(S,\mathbb{Q}(1))(-) \otimes \mathbb{C} \ar[r] & {\mathbb{Q}(-1)_{dR}} \otimes \mathbb{C} \ar[r] \ar@/^1pc/[l]^{s_{dR}} & 0}$$
and we must describe the two sections $s_B$ and $s_{dR}$. This is
inspired by the considerations in \cite{HaLNM}, 4.3.2, and \cite{HaMM}, I. Note that we neglect $\psi_f$. Along the rules we must find a form
$\omega_{top} \in \Hom_K(\Lambda^2(\mathfrak{g}/\mathfrak{k}),\mathcal{C}^{\infty}(B(\mathbb{Q})\backslash G(\mathbb{A}))),$ 
whose cohomology class $[\omega_{top}]$ generates $H^2(\partial S(\mathbb{C})^{BS},\mathbb{Q}(1)) \otimes \mathbb{C}$, and where the involution $F_{\infty}$ acts by $-1$. By Lemma \ref{LEMGEN} we get
$\omega_{top} := 2 \pi i \cdot  \delta x_1 \wedge \delta x_2,$
and we have to note that, because we have $\mathbb{Q}(1)$-coefficients, the involution $F_{\infty}$ acts indeed by $-1$. Then $[\Eis(2 \pi i \cdot \delta x_1 \wedge \delta x_2)]$ gives us $s_B(\mathbf{1}_{B})$.

Furthermore, we must find a form
$\omega_{hol} \in \Hom_K(\Lambda^2(\mathfrak{g}/\mathfrak{k}),\mathcal{C}^{\infty}(B(\mathbb{Q})\backslash G(\mathbb{A}))),$ 
whose cohomology class $[\omega_{hol}]$ generates $H^2(\partial S(\mathbb{C})^{BS},\mathbb{Q}(1)) \otimes \mathbb{C}$, and additionally it must be in $F^1H^2_{dR}(S(\mathbb{C}),\mathbb{Q}(1)) \otimes \mathbb{C}$ (cf. Section \ref{MMwC}). 
This is fulfilled by 
$[\omega_{hol}] =[ \delta z_1 \wedge \delta z_2],$ 
where we again have to put into account our Tate twist by $\mathbb{Q}(1)$. So actually, we look at $F^2H^2_{dR}(S(\mathbb{C}),\mathbb{Q}) \otimes \mathbb{C}$. And the Hodge filtration tells us that this is  $H^0(S(\mathbb{C}), \Omega^2).$
 Then $[\Eis (2 \pi i \cdot \delta z_1 \wedge \delta z_2)]$ gives us
 $s_{dR}(\mathbf{1}_{dR})$. Let us make this more precise. We start with $[2 \pi i \cdot \delta x_1 \wedge \delta x_2]$ as the generator $\mathbf{1}_{B} \in \mathbb{Q}(-1)_B \otimes \mathbb{C}$. This goes to 
$[\delta z_1 \wedge \delta z_2] = (2 \pi i)^{-1} \cdot \mathbf{1}_{dR} \in \mathbb{Q}(-1)_{dR} \otimes \mathbb{C}.$
And we have to compute $I_{\infty}^{-1}([ \delta z_1 \wedge \delta z_2]) \in H^2_{\CHE,B}(S,\mathbb{Q}(1))(-) \otimes \mathbb{C}$. For this we go around the square again in the other direction, i.e.
$I_{\infty}^{-1}([\delta z_1 \wedge \delta z_2]) = [\Eis (2 \pi i \cdot \delta z_1 \wedge \delta z_2)].$ 
So we get the difference 
$[\Eis (2 \pi i \cdot \delta z_1 \wedge \delta z_2)] - [\Eis (2 \pi i \cdot \delta x_1 \wedge \delta x_2)],$
which is in $\mathbb{Q}(0)\chi_{D,B} \otimes \mathbb{C}.$ 
And we are left with the multiplication with the Gaus\ss -sum $(\sqrt{D})^{-1}$, i.e.
$[H^2_{\CHE,\infty}(S,\mathbb{Q}(1))(-)]$ is $\frac{2 \pi i}{\sqrt{D}} \cdot
\left([\Eis ( \delta z_1 \wedge \delta z_2)] - [\Eis ( \delta x_1 \wedge
  \delta x_2)]\right).$ The first term $[\Eis(\delta z_1 \wedge \delta z_2)]$
is computed by Theorem \ref{MT2}. As $\delta z_1 \wedge \delta z_2 = \delta
x_1 \wedge \delta x_2 + i(\delta x_1 \wedge \delta y_2 - \delta x_2 \wedge
\delta y_1) + \delta y_1 \wedge \delta y_2,$ 
we have $\alpha_1 = i = - \alpha_2.$
Therefore 
$[\Eis(\delta z_1 \wedge \delta z_2)] = [\Eis( \delta x_1 \wedge \delta x_2)]
+ i \cdot  \frac{\sqrt{D} \cdot \log \varepsilon}{4 \pi  \cdot \zeta_F(-1)}
\left( c_1(L_1) - c_1(L_2)\right).$ 
Since the bottom $\mathbb{Q}(0)\chi_{D,B} \otimes \mathbb{C}$ is generated by
$2 \pi i \cdot (c_1(L_1) - c_1(L_2))$, the extension class is this multiple of $2 \pi i \cdot ( c_1(L_1) - c_1(L_2))$. So $[H^2_{\CHE,\infty}(S,\mathbb{Q}(1))(-)]  = \frac{i \cdot\log \varepsilon}{4
  \pi \cdot \zeta_F(-1)}.$ 
If we choose $\frac{1}{2 \pi i}$ as a basis for $i\mathbb{R}$, we are left with
$
[H^2_{\CHE,\infty}(S,\mathbb{Q}(1))(-)] = - \frac{\log \varepsilon}{2 \cdot \zeta_F(-1)}.
$
\end{proof}
Again the realisation is that of the Kummer motive $K \langle
\widetilde{\varepsilon} \rangle$, i.e. 
$K \langle \widetilde{\varepsilon} \rangle_{\infty} = [H^2_{\CHE,\infty}(S,\mathbb{Q}(1))(-)] = \log \widetilde{\varepsilon}.$

\section{Kummer-One-Motives} \label{K1-motive}
In this chapter we give even more evidence that the Kummer-Chern-Eisenstein
motive $[H^2_{\CHE}(S,\mathbb{Q}(1))(-)]$ is the Kummer motive $K
\langle \widetilde{\varepsilon} \rangle$ attached to $ \widetilde{\varepsilon}$. This relies on the
observation that $K \langle
\widetilde{\varepsilon} \rangle$ is isomorphic to the one-motive $M_{\widetilde{\varepsilon}}$
attached to $ \widetilde{\varepsilon}$ in the sense of \cite{De-HodgeIII}. 
I show how our  
$[H^2_{\CHE}(S,\mathbb{Q}(1))(-)]$
gives rise to the Kummer-$1$-motive $M_{\widetilde{\varepsilon}}$ attached to the
element $\widetilde{\varepsilon}$. There is a third (one-)motive related to our surface $S$, the Hodge-one-motive $\eta_{S}$. It corresponds to the Hodge structure of the cohomology of $S$. We meet this in Section \ref{SSHOM}.
\subsection{Kummer-Chern-Eisenstein vs. Kummer-One-Motives}\label{OM}%
Let us briefly recall the definition of a $1$-motive in the sense of Deligne
(see \cite{De-HodgeIII}, 10). Since we are in a very easy particular
situation, we do not need the general theory.

A one-motive (or $1$-motive) over a (algebraically closed) field $k$ is
defined by a complex $[X\stackrel{u}{\rightarrow}G]$, where $X$ is a
finitely generated free $\mathbb{Z}$-module, $G$ is a semi-abelian variety over
$k$, i.e. an extension of an abelian variety by a torus, and $u:X\rightarrow G(k)$ a group homomorphism. 
\begin{Rmk}\label{RMKOMMG}
If $k$ is not algebraically closed, but still a perfect field, one claims a $\text{Gal}(\overline{k}/k)$-action on $X$ and $G$, and the morphism $u$ is supposed to be morphism of $\text{Gal}(\overline{k}/k)$-modules. 
\end{Rmk}
Such a one-motive $M$ gives rise to a motive 
$T(M)=(T_B(M),T_{dR}(M),T_l(M)),$ see \cite{De-HodgeIII}, 10.1. 
For example: 
$T([\mathbb{Z}\rightarrow 0])=\mathbb{Z}(0)$, $T([0\rightarrow \mathbb{G}_m]) =
\mathbb{Z}(1)$ or $T([\mathbb{Z}\rightarrow \mathbb{G}_m]) \in
\Ext^1_{\mathcal{MM}_k}(\mathbb{Z}(0),\mathbb{Z}(1))$. We would like
 to see that the Kummer-Chern-Eisenstein motive in $\Ext^1_{\mathcal{MM}_{\mathbb{Q}}}(\mathbb{Q}(-1),\mathbb{Q}(0)\chi_D),$
 which was constructed in Chapter \ref{Construction}, actually comes from a one-motive, i.e. there is a $1$-motive $\mathcal{M}$ such that
$T(\mathcal{M}) \otimes \mathbb{Q} = [H^2_{\CHE}(S,\mathbb{Q}(1))(-)].$ 
Here we have to take into account the Galois action, that is given by the character $\chi_D$.

Each $1$-motive $[X\rightarrow G]$ is an extension in the category of $1$ -
motives $1\mathcal{M}_k$ i.e.
$0\rightarrow [0\rightarrow G]\rightarrow [X\rightarrow G]\rightarrow
[X\rightarrow 0]\rightarrow 0.$ 
On the other hand, each extension of $[\mathbb{Z}\rightarrow 0]$ by $[0\rightarrow \mathbb{G}_m]$, i.e. an element in  $\Ext^1_{1\mathcal{M}_k}([\mathbb{Z}\rightarrow 0],[0\rightarrow \mathbb{G}_m])$, is of the form $M=[\mathbb{Z}\stackrel{u}{\rightarrow} \mathbb{G}_m]$.
\begin{Def}
We call an element in $\Ext^1_{1\mathcal{M}_k}([\mathbb{Z}\rightarrow
0],[0\rightarrow \mathbb{G}_m])$ a Kummer-$1$-motive and denote it by $M_t =
[\mathbb{Z} \stackrel{u}{\rightarrow} \mathbb{G}_m],$ that is $u(1)=t$.
\end{Def}
Here we know that such an $M_t$ is uniquely determined by $t\in \mathbb{G}_m(k)$, i.e.
we have 
$\Ext^1_{1\mathcal{M}_k}([\mathbb{Z}\rightarrow 0],[0\rightarrow \mathbb{G}_m]) = k^*.$
Recall (Section \ref{MMwC}) that we have for the category $\mathcal{MM}_k$ of mixed motives over $k$, only a conjecture of such an identification, i.e. we can only expect
$$ \Ext^1_{\mathcal{MM}_k}(T([\mathbb{Z} \rightarrow 0]) \otimes \mathbb{Q},T([0\rightarrow \mathbb{G}_m]) \otimes \mathbb{Q}) \subset \Ext^1_{\mathcal{MM}_k}(\mathbb{Q}(0),\mathbb{Q}(1)).$$
In other words, the Kummer-one-motives and the Kummer motives should form the same subcategory, and indeed we have
\begin{Lem}\label{LEMK=K}

Let $a \in k^*$. Let $K \langle a \rangle$ the Kummer motive attached to $a$, and $M_a$ the Kummer-one-motive attached to $a$. Then $K \langle a \rangle = M_a.$
\end{Lem}
\begin{proof}
The construction of $K \langle a \rangle$ is explained in detail in \cite{TK},
3.1 (we sketched this construction in the beginning of Section
\ref{MMwC}). But this is exactly the same construction for $M_a$ as in
\cite{De-HodgeIII}, 10.3. Compare also \cite{AJS}, 2.7.
\end{proof}
This lemma puts us in an even better situation to conclude that the our
Kummer-Chern-Eisenstein motive is indeed a Kummer motive. So we must find the one
dimensional $\mathbb{Z}$-module $X$, the multiplicative group $\mathbb{G}_m$
and the map $u$. Let us look at the dual situation (Section \ref{dual
  motive}). The dual $[H^2_{\CHE}(S,\mathbb{Q}(1))(-)]^{\vee}$ is in
$\Ext^1_{\mathcal{MM}_{\mathbb{Q}}}(\mathbb{Q}(0)\chi_D,\mathbb{Q}(1))$ and by the above third example we know
$$ \Ext^1_{\mathcal{MM}_{\mathbb{Q}}}(T([\mathbb{Z}(\chi_D)\rightarrow 0]) \otimes \mathbb{Q},T([0\rightarrow \mathbb{G}_m]) \otimes \mathbb{Q}) \subset \Ext^1_{\mathcal{MM}_{\mathbb{Q}}}(\mathbb{Q}(0)\chi_D,\mathbb{Q}(1)).$$
Again we must be aware of the action of Galois given by the character $\chi_D$, i.e. the realisation $T([\mathbb{Z}(\chi_D)\rightarrow 0]) \otimes \mathbb{Q}$ is the Dirichlet motive $\mathbb{Q}(0)\chi_D$.

Recall (by Section \ref{Construction} and \ref{dual motive}) that the top of
the extension $\mathbb{Q}(0)\chi_D$ is generated up to a constant by the first Chern class 
$c_1(L) = c_1(L_1^{-1} \otimes L_2) \in H^2_{!}(S,\mathbb{Q}(1)),$ 
see Remark \ref{RMKDM}. The middle
$H^2_{\CHE}(S,\mathbb{Q}(1))(-)^{\vee}$ sits in $H^2_c(S,\mathbb{Q}(1))$ - see loc. cit. 
The bottom $\mathbb{Q}(1)$ comes from the cohomology group $H^1(\widetilde{S},\mathbf{R}j_*\mathbb{Q}/j_!\mathbb{Q})$ twisted by $\mathbb{Q}(1)$.

Now consider the restriction map 
$u: \text{Pic}(\widetilde{S}) \rightarrow \text{Pic}(\widetilde{S}_{\infty}).$
In $\text{Pic}(\widetilde{S})$ we have the element
$\widetilde{L}=\widetilde{L_1}^{-1} \otimes \widetilde{L_2}$, whose Chern
class $c_1(\widetilde{L_1}^{-1} \otimes \widetilde{L_2})$ generates
$\mathbb{Q}(0)\chi_D$ (Lemma \ref{MOTKK}). Furthermore, we know that this has got trivial Chern class on the boundary, i.e. 
$u(\widetilde{L_1}^{-1} \otimes \widetilde{L_2}) \in \text{Pic}^0(\widetilde{S}_{\infty}),$ 
and finally we know that
$\text{Pic}^0(\widetilde{S}_{\infty}) \simeq \mathbb{G}_m,$ that is we get a $1$-motive 
$[\mathbb{Z}(\chi_D) \cdot \widetilde{L} \stackrel{u}{\rightarrow} \text{Pic}^0(\widetilde{S}_{\infty})].$
\begin{Lem}\label{LME}
Let $\widetilde{L}=\widetilde{L_1}^{-1} \otimes \widetilde{L_2} \in \text{Pic}(\widetilde{S})$ and $u : \text{Pic}(\widetilde{S}) \rightarrow \text{Pic}(\widetilde{S}_{\infty})$ be as above. Let $\varepsilon \in \mathcal{O}_F^*$ be as fixed in the very beginning. Then the $1$-motive $[\mathbb{Z}(\chi_D) \cdot \widetilde{L} \stackrel{u}{\rightarrow} \text{Pic}^0(\widetilde{S}_{\infty})]$ is equal to the Kummer-$1$-motive attached to $\varepsilon^{-2}$.
\end{Lem} 
\begin{proof}
In Lemma \ref{LEMCLB} (see also Corollary \ref{i*e}) we proved that $\widetilde{L_1}^{-1} \otimes \widetilde{L_2}$ goes via the restriction map to $\varepsilon^{-2}$.
\end{proof}
So write  
$M_{\varepsilon^{-2}} = [\mathbb{Z}(\chi_D) \cdot \widetilde{L} \stackrel{u}{\rightarrow} \text{Pic}^0(\widetilde{S}_{\infty})]$
as an extension
$$0 \rightarrow [0 \rightarrow  \text{Pic}^0(\widetilde{S}_{\infty})] \rightarrow M_{\varepsilon^{-2}} \rightarrow [\mathbb{Z}(\chi_D) \cdot \widetilde{L} \rightarrow 0] \rightarrow 0,$$
i.e. as an element in $\Ext^1_{1\mathcal{M}_{\mathbb{Q}}}([\mathbb{Z}(\chi_D)\rightarrow 0],[0\rightarrow \mathbb{G}_m])$.
\begin{Lem}\label{1u1}
Let $ \text{Pic}(\widetilde{S})$, $\text{Pic}^0(\widetilde{S}_{\infty})$ and
the line bundle $\widetilde{L} \in \text{Pic}(\widetilde{S})$ as above. Let
$\mathbb{Q}(0)\chi_D$ and
$H^1(\widetilde{S},\mathbf{R}j_*\mathbb{Q}/j_!\mathbb{Q})$ be the motives as
in Section \ref{Construction}. Then $T([\mathbb{Z}(\chi_D) \cdot
\widetilde{L} \rightarrow 0]) \otimes \mathbb{Q}$ is isomorphic to
$\mathbb{Q}(0)\chi_D$. And there is an isomorphism
$T([0 \rightarrow \text{Pic}^0(\widetilde{S}_{\infty})]) \otimes \mathbb{Q} \simeq H^1(\widetilde{S},\mathbf{R}j_*\mathbb{Q}/j_!\mathbb{Q}) \otimes  \mathbb{Q}(1).$
\end{Lem}
\begin{proof}
We know that $\text{Pic}^0(\widetilde{S})$ is trivial, i.e. the Chern class
map $c_1$ sends the
line bundle $\widetilde{L}$ uniquely to $c_1(\widetilde{L}) \in H^2(\widetilde{S},\mathbb{Q}(1))$. Since this class
generates $\mathbb{Q}(0)\chi_D$, we get that $T([\mathbb{Z}(\chi_D) \cdot
\widetilde{L} \rightarrow 0]) \otimes \mathbb{Q}$ is isomorphic to
$\mathbb{Q}(0)\chi_D$. We know by definition that 
$T([0 \rightarrow \text{Pic}^0(\widetilde{S}_{\infty})]) \otimes \mathbb{Q} =
T([0 \rightarrow \mathbb{G}_m]) \otimes \mathbb{Q} = \mathbb{Q}(1).$ We observe that our $\text{Pic}^0(\widetilde{S}_{\infty})$ is in $H^1(\widetilde{S}_{\infty},\mathcal{O}_{\widetilde{S}_{\infty}}^*)$, i.e. 
$T([0 \rightarrow \text{Pic}^0(\widetilde{S}_{\infty})]) \otimes \mathbb{Q} \simeq  H^1(\widetilde{S}_{\infty},\mathbb{Q}(1)).$
On the other hand, we have (Lemma \ref{MOTKK}) 
$H^1(\widetilde{S}_{\infty},\mathbb{Q}(1)) = H^1(\widetilde{S},\mathbf{R}j_*\mathbb{Q}/j_!\mathbb{Q}) \otimes \mathbb{Q}(1) =  \mathbb{Q}(1),$
i.e. 
$T([0 \rightarrow \text{Pic}^0(\widetilde{S}_{\infty})]) \otimes \mathbb{Q} \simeq  H^1(\widetilde{S},\mathbf{R}j_*\mathbb{Q}/j_!\mathbb{Q}) \otimes \mathbb{Q}(1).$
\end{proof}
Now we conclude that the motive
$T([\mathbb{Z}(\chi_D) \cdot \widetilde{L} \stackrel{u}{\rightarrow} \text{Pic}^0(\widetilde{S}_{\infty})]) \otimes \mathbb{Q}$ is isomorphic to $0 \rightarrow \mathbb{Q}(1) \rightarrow T(M_{\varepsilon^{-2}}) \otimes \mathbb{Q} \rightarrow \mathbb{Q}(0)\chi_D \rightarrow 0,$
i.e. an element in
$\Ext^1_{\mathcal{MM}_{\mathbb{Q}}}(\mathbb{Q}(0)\chi_D,\mathbb{Q}(1)).$ 
And if we normalise (Remark \ref{RMKDM}) the dual one is 
$0 \rightarrow \mathbb{Q}(0)\chi_D \rightarrow
T(M_{\widetilde{\varepsilon}})\otimes \mathbb{Q} \rightarrow \mathbb{Q}(-1)
\rightarrow 0$. 
Now I would like to claim even more that indeed 
$[H^2_{\CHE}(S,\mathbb{Q}(1))(-)] \simeq T(M_{\widetilde{\varepsilon}}) \otimes
\mathbb{Q}.$ We have this at present for the $l$-adic realisations.
\begin{Thm}\label{KCEMOM}
Let $[H^2_{\CHE}(S,\mathbb{Q}(1))(-)]$ be our Kummer-Chern-Eisenstein motive
and $T([\mathbb{Z}(\chi_D) \cdot \widetilde{L}
\stackrel{u}{\rightarrow} \text{Pic}^0(\widetilde{S}_{\infty})]^{\vee})
\otimes \mathbb{Q}$. Then the $l$-adic realisations are isomorphic 
$[H^2_{\CHE,l}(S,\mathbb{Q}(1))(-)] \simeq T_l([\mathbb{Z}(\chi_D) \cdot
\widetilde{L} \stackrel{u}{\rightarrow}
\text{Pic}^0(\widetilde{S}_{\infty})]^{\vee}) \otimes \mathbb{Q}_l.$ 
In particular, we have
$[H^2_{\CHE,l}(S,\mathbb{Q}(1))(-)] \simeq T_l(M_{\widetilde{\varepsilon}}) \otimes \mathbb{Q}_l.$
\end{Thm}
\begin{proof}
We add the map $H^1_{\acute{e} t}(i^*)$ on the bottom of the big diagram in Lemma \ref{BD} - see beginning of
Section \ref{liftings}. Then we find 
the $l$-adic realisations of the $1$-motive 
$T_l(M_{\varepsilon^{-2}}) \otimes \mathbb{Q}_l = T_l([\mathbb{Z}(\chi_D) \cdot \widetilde{L} \stackrel{u}{\rightarrow} \text{Pic}^0(\widetilde{S}_{\infty})]) \otimes \mathbb{Q}_l$ 
as 
$$\xymatrix{[H^1_{\acute{e}
    t}(\widetilde{S}\times\overline{\mathbb{Q}},\mathbb{G}_{m}) \ar[r]
  \ar@/_2pc/@{-->}[rr]_{u: \widetilde{L} \mapsto \varepsilon^{-2}} &
  H^1_{\acute{e}
    t}(\widetilde{S}_{\infty}\times\overline{\mathbb{Q}},i^*\mathbb{G}_{m,\widetilde{S}})
    \ar[r] & H^1_{\acute{e}
    t}(\widetilde{S}_{\infty}\times\overline{\mathbb{Q}},\mathbb{G}_{m,\widetilde{S}_{\infty}})]
    }.$$ 
We know that the dual motive sits in the following sequence (see Lemma \ref{MOTKK})
$$0  \rightarrow {
  \underset{n}{\varprojlim}H^1_{\acute{e}t}(\widetilde{S}_{\infty}\times\overline{\mathbb{Q}},\boldsymbol{\mu}_{l^n})
  \otimes \mathbb{Q}_l } \rightarrow  { \underset{n}{\varprojlim}
  H^2_{\acute{e} t,c}(S\times\overline{\mathbb{Q}},\boldsymbol{\mu}_{l^n})
  \otimes \mathbb{Q}_l} \rightarrow  {\text{Im}(f_1)} \rightarrow 0,$$ where
  this sequence itself is part of the big diagram of Lemma \ref{BD}.  
Recall furthermore, that the extension class of
$[H^2_{\CHE,l}(S,\mathbb{Q}(1))(-)]^{\vee}$ is $\varepsilon^{-2}$ and is in the bottom $\mathbb{Q}_l(1)$. By Lemma \ref{1u1} we see that the generator $\widetilde{L}$ of  $T_l([\mathbb{Z}(\chi_D) \cdot \widetilde{L} \rightarrow 0]) \otimes \mathbb{Q}_l$ goes via the Chern class map uniquely to $c_1(\widetilde{L}) \in {\text{Im}(f_1)}$, and this class generates ${\mathbb{Q}_l(0) \chi_D}$.
 On the other hand, we know (Lemma \ref{LME}) that this generator $\widetilde{L}$ maps under $u = H^1_{\acute{e} t}(i^*) \circ R$ to $T_l([0 \rightarrow \text{Pic}^0(\widetilde{S}_{\infty})]) \otimes \mathbb{Q}_l$, and its image is  
$u(\widetilde{L}) = (H^1_{\acute{e} t}(i^*) \circ R)(\widetilde{L}) = (H^1_{\acute{e} t}(i^*) \circ \varrho)(c_1(\widetilde{L})) = \varepsilon^{-2}.$ 
Now the $l$-adic realisation  $T_l([0 \rightarrow \text{Pic}^0(\widetilde{S}_{\infty})]) \otimes \mathbb{Q}_l$ is the Tate module of $\text{Pic}^0(\widetilde{S}_{\infty})$, i.e.
$\mathfrak{T}_l( \text{Pic}^0(\widetilde{S}_{\infty})) = \mathfrak{T}_l(\mathbb{G}_m),$ 
but this is exactly 
$\underset{n}{\varprojlim} H^1_{\acute{e} t}(\widetilde{S}_{\infty} \times \overline{\mathbb{Q}},\boldsymbol{\mu}_{l^n}) \otimes \mathbb{Q}_l,$ 
therefore the $l$-adic realisation class of the one-motive is in the bottom
$\mathbb{Q}_l(1)$. Now we use the same diagram chases as in the proof of Theorem \ref{MTladic}. This gives us 
$ [H^2_{\CHE,l}(S,\mathbb{Q}(1))(-)]^{\vee}  \simeq T_l([\mathbb{Z}(\chi_D) \cdot \widetilde{L} \stackrel{u}{\rightarrow} \text{Pic}^0(\widetilde{S}_{\infty})]^{\vee})  \otimes \mathbb{Q}_l  =   T_l(M_{\varepsilon^{-2}}) \otimes \mathbb{Q}_l.$
And we are left with dualising.
\end{proof}
To get the deeper result that $[H^2_{\CHE}(S,\mathbb{Q}(1))(-)] \simeq
T(M_{\widetilde{\varepsilon}}) \otimes \mathbb{Q}$, we have to assume that
$[H^2_{\CHE}(S,\mathbb{Q}(1))(-)]$ is indeed a Kummer motive. Then we can
refer to \cite{Ja}, Theorem 4.3, which allows that it is sufficient to look at the $l$-adic realisations. This theorem says that two one-motives are isomorphic, if and only if the $l$-adic realisations $T_l$ are isomorphic. So if  $[H^2_{\CHE}(S,\mathbb{Q}(1))(-)]$ is a Kummer motive, then it has to be $M_{\widetilde{\varepsilon}}$.
\subsection{The Hodge-One-Motive of an Algebraic Surface}\label{SSHOM}%
In this final section we show how our Kummer-Chern-Eisenstein
motive and the $1$-motive
$M_{\widetilde{\varepsilon}}$ fit into the picture of \cite{Ca}. In
loc. cit. one considers the case of a complex surface. Since
$[H^2_{\CHE}(S,\mathbb{Q}(1))(-)]$ is defined over $\mathbb{Q}$, we need a
generalisation to arbitrary base fields. 
This has been done (even for higher dimensions) in the works of L. Barbieri-Viale, et.al. \cite{BVRS} and independently of N. Ramachandran \cite{Ram}. These tell us that the Hodge-one-motive is indeed defined over $\mathbb{Q}$. Hence in the following we refer to \cite{Ca}, but keep the others in mind. If one wants to avoid these generalisations, one can just look at the complex situation.

The starting point is Deligne's observation (\cite{De-HodgeIII}, 10.1.3) that
the category of $1$-motives is equivalent to the category of torsion-free
mixed Hodge structures of length one. Consider the largest Hodge substructure of $H^2(X,\mathbb{Q}(1))$, which is of
type $\{(0,0),(-1,-1),(0,-1),(-1,0)\},$ where $X$ is a complex algebraic variety. Then, by the above equivalence of categories, there is a unique $1$-motive $\eta_X$ corresponding to this Hodge structure, which is called the Hodge-$1$-motive.
Now \cite{Ca}, Theorem K, delivers that for a complex algebraic surface $S$
there is a geometric construction of $\eta_S$, i.e. it is isomorphic to a $1$
- motive, called the \emph{trace-1-motive} $\tau_S$ of $S$. Here in our case
of the Hilbert modular surface $S$ it comes down to an easy situation, as we
described in the last Section \ref{OM} (compare also \cite{Ca}, Chapter 15). 

Let
$\text{NS}(\widetilde{S},\widetilde{S}_{\infty}) \subset \text{NS}(\widetilde{S})$
be the subgroup of the Neron-Severi group $\text{NS}(\widetilde{S})$,
consisting of those cycles on the compact surface $\widetilde{S}$, which have
trivial Chern class on the boundary divisor $\widetilde{S}_{\infty}$. Note
that this is a finitely generated $\mathbb{Z}$-module, and furthermore,
since $\text{Pic}^0(\widetilde{S})$ vanishes, we have that
$\text{Pic}(\widetilde{S}) = \text{NS}(\widetilde{S})$. Then we get, by
\cite{Ca}, Theorem K, the Hodge-one -motive of $H^2(\widetilde{S} - \widetilde{S}_{\infty},\mathbb{Q}) = H^2(S,\mathbb{Q})$ via the restriction map
$\eta_S = \tau_{S}:\text{NS}(\widetilde{S},\widetilde{S}_{\infty})\rightarrow
\text{Pic}^0(\widetilde{S}_{\infty}).$ On the other hand, we have got our Kummer-$1$-motive
$[\mathbb{Z}(\chi_D) \cdot \widetilde{L} \stackrel{u}{\rightarrow} \text{Pic}^0(\widetilde{S}_{\infty})] \simeq M_{\varepsilon^{-2}}.$
We can consider $M_{\varepsilon^{-2}}$ as a submotive of $\tau_{S} = \eta_S$, since the generator $\widetilde{L}$ is mapped uniquely to $c_1(\widetilde{L}) \in \text{NS}(\widetilde{S},\widetilde{S}_{\infty})$ (compare Lemma \ref{1u1}). Now according to Theorem \ref{KCEMOM} above this should be the dual of our Kummer-Chern-Eisenstein motive $[H^2_{\CHE}(S,\mathbb{Q}(1))(-)]^{\vee}$. To get the things in order, one must look briefly at the weight filtration of
the Hodge structures. First we observe, since we deal with $\mathbb{Q}(1)$ -
coefficients, that we have to consider a substructure of type
$\{(0,0),(-1,-1),(0,-1),(-1,0)\}$ for the Hodge-$1$-motive $\eta_S$. The
weight filtration $\{W_{\bullet} H^2(S,\mathbb{Q}(1))\}$ of
$H^2(S,\mathbb{Q}(1))$ with $\mathbb{Q}(1)$-coefficients is (see \cite{vdG}, VI.1)
$$W_k H^2(S,\mathbb{Q}(1)) = 
\begin{cases}
0 & k \leq 0 \\
H^2_!(S,\mathbb{Q}(1)) & k = 1 = 2 \\
0 & k = 3 \\
H^2(S,\mathbb{Q}(1)) & k \geq  4.
\end{cases}$$ 
Therefore we must focus on the $\text{Gr}_1^W$-part, and this is just
$H^2_!(S,\mathbb{Q}(1))$. For this we have the Hodge filtration $\{F^{\bullet}
H^2_!(S,\mathbb{Q}(1)) \otimes \mathbb{C}\}$ (compare e.g. loc. cit. Proposition 1.2). And in view of the one-motives, we must take the
$(0,0)$-part (note our Tate twist). It consists of the two Chern classes
$c_1(L_1)$ and $c_1(L_2)$ and some $(1,1)$-forms coming from the cuspidal
part. In particular, we find there our motive $\mathbb{Q}(0) \chi_D$. One should think of a picture like this for $\text{Gr}_1^W$ 
$$\xymatrix{1 \ar@{-}[d] \ar@{--} '[dr]^{F^{\bullet}} '[ddrr] \ar@{-}[rr] &  & (1,1) \ar@{-}[dd] \\0 \ar@{-}[d] \ar@{-}[r] \ar@{}[dr] |{\eta_{S}} & *+[F--]{\mathbb{Q}(0) \chi_D} \ar@{-}[d] &  \\ -1 \ar@{-}[r] & 0 \ar@{-}[r] & 1. }$$
Moreover, we cannot expect to get the whole Hodge-one-motive.
\begin{Thm}\label{THMKOM}

Let $S$ be our Hilbert modular surface. Consider the Kummer-Chern-Eisenstein motive $[H^2_{\CHE}(S,\mathbb{Q}(1))(-)]$,
the Kummer-(one)-motive $M_{\varepsilon^{-2}}$ attached to the element
$\varepsilon^{-2}$ and
the Hodge-one-motive $\eta_{S}$ attached to the Hodge structure of
$H^2(S,\mathbb{Q}(1))$. Then the Kummer-$1$-motive attached to
$\varepsilon^{-2}$ is isomorphic to a submotive of the Hodge-one-motive
$\eta_{S}$. In particular, the dual of the Kummer-Chern-Eisenstein motive $[H^2_{\CHE}(S,\mathbb{Q}(1))(-)]^{\vee}$ is isomorphic to a submotive of the realisation $T(\eta_S) \otimes \mathbb{Q}$ of the Hodge-one-motive.
\end{Thm}
\begin{proof}
As described above we know by \cite{Ca}, Theorem K, that the Hodge-one-motive $\eta_{S}$ is isomorphic to,
$\tau_{S}:\text{NS}(\widetilde{S},\widetilde{S}_{\infty})\rightarrow \text{Pic}^0(\widetilde{S}_{\infty}),$
i.e. by Lemma \ref{LME} we get the submotive $[\mathbb{Z}(\chi_D) \cdot \widetilde{L} \stackrel{u}{\rightarrow} \text{Pic}^0(\widetilde{S}_{\infty})] \simeq M_{\varepsilon^{-2}}$.
\end{proof}
 ALEXANDER CASPAR \\
ETH Z\"urich - Departement Mathematik \\
R\"amistrasse 101, CH-8092 Z\"urich \\
Switzerland\\
caspar@math.ethz.ch
\addcontentsline{toc}{section}{References}  


\end{document}